\documentclass[twoside,parskip]{article}
\pdfoutput=1

	\usepackage{amsmath}				
	\usepackage{amscd}					
	\usepackage{amsthm,thmtools}        
	\usepackage{amssymb}				
	
	\usepackage[titletoc]{appendix}		
	\usepackage[
		style=alphabetic,				
		sorting=nyt,					
		maxnames=5,						
		minnames=3,						
		backend=bibtex8,				
		firstinits=true,				
		doi=false,						
	]{biblatex}
	\usepackage{bm}						
	\usepackage{color}					
	\usepackage{enumitem}				
	\usepackage{etoolbox}				
	\usepackage{fancyhdr}				
	\usepackage{float}					
	\usepackage[T1]{fontenc}			
	\usepackage{gensymb}				
	\usepackage[a4paper,
		outer=2cm,
		inner=2cm,
		top=2.5cm,
		bottom=2cm]{geometry}
	\usepackage[utf8]{inputenc}			
	\usepackage{ifmtarg}				
	\usepackage{mathrsfs}				
	\usepackage{makeidx}				
	\usepackage{multicol}				
	\usepackage{parskip}				
	\usepackage[refpage,intoc]{nomencl}		
	\usepackage{slashed}				
	\usepackage{stmaryrd}				
	\usepackage{tensor}					
	\usepackage{textcomp}				
	\usepackage{titlesec}               
	\usepackage{todonotes}
	\usepackage{tikz}                   
	\usepackage{tocloft}				
	\usepackage{url}					
	\usepackage[all,cmtip]{xy}			
	
	\usepackage{csquotes}				

	\usepackage[breaklinks,pageanchor]{hyperref}		
	\usepackage[capitalize]{cleveref}	

	\bibliography{refs} 
	\DeclareFieldFormat{postnote}{#1}
	\DeclareFieldFormat{multipostnote}{#1}
	
		\numberwithin{equation}{section}
		\crefformat{equation}{(#2#1#3)} 

	\fancypagestyle{plain}{%
		\fancyhf{}%
		\fancyfoot[C]{}%
	}

		\let\oldtitle\title
		\renewcommand{\title}[1]{\oldtitle{#1}\newcommand{\mythetitle}{#1}}
		\let\oldauthor\author
		\renewcommand{\author}[1]{\oldauthor{#1}\newcommand{\mytheauthor}{#1}}	
		\definecolor{linkcolor}{RGB}{00,10,138} 
		\hypersetup{
			colorlinks=true,
			breaklinks=true,
			linkcolor=linkcolor,
			menucolor=linkcolor,
			citecolor=linkcolor,
			filecolor=linkcolor,
			urlcolor=linkcolor,
			frenchlinks=false
		}
	
		\renewcommand{\nomname}{List of Symbols}
		\patchcmd{\thenomenclature}
		  {\chapter*{\nomname}}
		  {\phantomsection\chapter*{\nomname}\label{SectNomenclature}}
		  {}
		  {}

		\titleformat{\section}[block]{\large\bfseries\scshape\filcenter\textsc}{\thesection. }{0pt}{}
		\titleformat{\subsection}{\bfseries}{\thesubsection. }{0pt}{}

	\declaretheoremstyle[
	spaceabove=6pt, spacebelow=6pt,
	headfont=\normalfont\bfseries,
	notefont=\bfseries, notebraces={(}{)},
	bodyfont=\normalfont,
	postheadspace=1em,
	qed=$\lozenge$
	]{mystylesingle}

	\declaretheorem[style=mystylesingle,name=Main Theorem]{MainThm}

	\declaretheoremstyle[
	spaceabove=6pt, spacebelow=6pt,
	headfont=\normalfont\bfseries,
	notefont=\bfseries, notebraces={(}{)},
	bodyfont=\normalfont,
	postheadspace=1em,
	qed=$\lozenge$,
	]{mythmstyle}

	\declaretheoremstyle[
	spaceabove=6pt, spacebelow=6pt,
	headfont=\normalfont\bfseries,
	notefont=\mdseries, notebraces={(}{)},
	bodyfont=\normalfont,
	postheadspace=1em,
	qed=$\square$,
	]{myprfstyle}

	\declaretheorem[style=mythmstyle,name=Definition,numberwithin=section]{Def}
	
	\declaretheorem[style=mythmstyle,sibling=Def,name=Theorem]{Thm}

	\declaretheorem[style=mythmstyle,sibling=Def,name=Lemma]{Lem}
	
	\declaretheorem[style=mythmstyle,sibling=Def,name=Remark]{Rem}

	\declaretheorem[style=myprfstyle,numbered=no,name=Proof]{Prf}

	\crefname{Def}{Definition}{Definitions}
	\crefname{Exm}{Example}{Examples}
	\crefname{Con}{Convention}{Conventions}
	\crefname{Cor}{Corollary}{Corollaries}
	\crefname{Fct}{Fact}{Facts}
	\crefname{Lem}{Lemma}{Lemmas}
	\crefname{NumText}{Paragraph}{Paragraphs}
	\crefname{Rem}{Remark}{Remarks}
	\crefname{Cnj}{Conjecture}{Conjectures}
	\crefname{Exc}{Exercise}{Exercises}
	\crefname{MainThm}{Main Theorem}{Main Theorems}
	\crefname{Prb}{Problem}{Problems}
	\crefname{Qst}{Question}{Questions}
	\crefname{Thm}{Theorem}{Theorems}
	\crefname{Not}{Notation}{Notation}
	\crefname{Prf}{Proof}{Proofs}

		\newcommand{\fun}[1]{\bm{\mathrm{#1}}}
		
		\sbox0{$\mathbf{\xdef\mybffam{\the\fam}}$}
		\newcommand{\cat}[1]{\begingroup\fam\mybffam#1\endgroup}

		\newlist{category}{enumerate}{10}
		\setlist[category]{
			itemsep=1ex,
			topsep=0.5ex,
			parsep=0pt,
			leftmargin=3em,
			itemindent=5em,
			labelwidth=5em,
			labelsep=1em,
			align=categorylabel,
			label*=.\arabic*
		}
		\setlist[category,1]{label=\arabic*}
		\SetLabelAlign{categorylabel}{\bfseries{\scshape{#1}}:\hfil}	
	
		\newcommand{\DefMap}[4]{
			\begin{align*}
				\begin{array}{rcl}
					#1 & \to & #2 \\
					#3 & \mapsto & #4 
				\end{array} 
			\end{align*}
		}

	\newcommand{\emf}[1]{\textbf{\textit{#1}}}

	\newlist{caselist}{enumerate}{10}
	\setlist[caselist]{itemsep=1ex,topsep=1ex,parsep=0pt,leftmargin=0pt,%
	  labelwidth=*,labelsep=1ex,align=caselabel,label*=.\arabic*}
	\setlist[caselist,1]{label=\arabic*}

	\newcommand\thecaselabeltext{}
	\SetLabelAlign{caselabel}{\textsc{Case~#1}\thecaselabeltext:\hfil}

	\crefname{caselisti}{case}{cases}
	\loop\ifnum\count255<10
	   \advance\count255 1
	   \crefalias{caselist\romannumeral\count255}{caselisti}
	\repeat

	\makeatletter
	\def\ifempty#1{\def\@x{#1}\ifx\@x\@empty}
	\makeatother

	\newlist{steplist}{enumerate}{10}
	\setlist[steplist]{nolistsep,itemsep=0ex,topsep=0ex,partopsep=0px,parsep=0pt,leftmargin=0pt,%
	  labelwidth=*,labelsep=1ex,align=steplabel,label*=.\arabic*}
	\setlist[steplist,1]{label=\arabic*}

	\newcommand\thesteplabeltext{}
	\SetLabelAlign{steplabel}{\textsc{Step~#1}\thesteplabeltext:\hfil}

	\crefname{steplisti}{Step}{steplist}
	\loop\ifnum\count255<10
	   \advance\count255 1
	   \crefalias{steplist\romannumeral\count255}{steplisti}
	\repeat

	\makeatletter
	\def\ifempty#1{\def\@x{#1}\ifx\@x\@empty}
	\makeatother

	\DeclareMathOperator{\Asym}{Asym} 				
	\DeclareMathOperator{\bbullet}{\overline{\bullet}}						
	\DeclareMathOperator{\Dirac}{\slashed{D}}		
	\DeclareMathOperator{\EDM}{EDM} 				
	\newcommand{\GLp}{\operatorname{GL}^{+}}		
	\newcommand{\GLtp}{\widetilde{\operatorname{GL}}^{+}}	
	\DeclareMathOperator{\GL}{GL} 					
	\DeclareMathOperator{\Iso}{Iso} 				
	\DeclareMathOperator{\Lor}{Lor} 				
	\DeclareMathOperator{\N}{\mathbb{N}}			
	\DeclareMathOperator{\Ric}{\operatorname{Ric}}
	\DeclareMathOperator{\R}{\mathbb{R}}			
	\DeclareMathOperator{\SO}{SO} 					
	\DeclareMathOperator{\Spin}{Spin} 				
	\DeclareMathOperator{\Sym}{Sym} 				
	\DeclareMathOperator{\id}{id} 					
	\DeclareMathOperator{\loc}{loc} 				
	\DeclareMathOperator{\scal}{scal}				
	\DeclareMathOperator{\tr}{tr} 					
	
	\renewcommand{\Re}{\operatorname{Re}} 					
	\renewcommand{\S}{\operatorname{S}}					

	\title{A Universal Spinor Bundle and the Einstein-Dirac-Maxwell Equation as a Variational Theory}
	\author{Olaf Müller, Nikolai Nowaczyk}
	
	\makeindex
	\usepackage[totoc]{idxlayout}
	
	\makenomenclature

\begin{document}
	\thispagestyle{empty}
		\begin{center}
			{\Large \bfseries \mythetitle} \\
			\vspace{0.3cm}
			\mytheauthor  \\
			\vspace{0.3cm}
			\today \\
		\end{center}
	
		\textbf{Abstract.} 
		Not only the Dirac operator, but also the spinor bundle of a pseudo-Riemannian manifold depends on the underlying metric. This leads to technical difficulties in the study of problems where many metrics are involved, for instance in variational theory. We construct a natural finite dimensional bundle, from which all the metric spinor bundles can be recovered including their extra structure. In the Lorentzian case, we also give some applications to Einstein-Dirac-Maxwell theory as a variational theory and show how to coherently define a maximal Cauchy development for this theory.
		
		\textbf{Keywords.} spin geometry, spinor bundle, jet spaces, natural constructions, Einstein-Dirac-Maxwell equation, Cauchy development
		
		\textbf{Mathematics Subject Classification 2010.} 35L03, 53C27, 58J45, 83C22
		
		
		\section{Introduction}

Spinor bundles are an important tool in differential geometry as well as in mathematical physics, where they model fermionic particles. As evident from their construction, spinor bundles depend on the underlying metric. This makes it difficult to compare spinor fields in spinor bundles formed with respect to different metrics. The same problem occurs if one wants to compare the Dirac operators, since their domains of definition are not the same. A way out of this dilemma is provided by systematically constructing identification isomorphisms, see \cite{BourgGaud,BaerGaudMor}. Although using these identification isomorphisms is sufficient for many applications, one might wonder if and how it is possible to define a natural finite dimensional bundle that is independent of the metric, but nevertheless captures the features of the spinor bundles formed with respect to all the various metrics. 

In the present article, we propose the following answer to this question: Let $(M, \Theta)$ be a smooth spin manifold of dimension $m$ and $\Theta:\GLtp M \to \GLp M$ be a topological spin structure on $M$. Let $\tau_M: TM \to M$ be the tangent bundle of $m$ and for any $r+s=m$, let $\pi^{r,s}: \S_{r,s}M \to M$ be the subbundle of $\tau^*_{M} \otimes \tau^*_M: T^*M \otimes T^* M \to M$, which over each point $x \in M$ consists of all non-degenerate symmetric bilinear forms $g_x$ on $T_x M$ of signature $(r,s)$. Here, $r$ denotes the dimension of a maximal positive definite subspace and $s$ of a maximal negative definite subspace. We define $\mathcal{S}_{r,s}(M) := \Gamma(\pi^{r,s})$ to be the space of all metrics on $M$ of signature $(r,s)$. For any metric $g \in \mathcal{S}_{r,s}(M)$, we obtain a \emf{metric spinor bundle}
\begin{align}
	\label{EqDefSpinorBundleClassic}
	\pi^g_M: \Sigma^g M := \Spin^g M \times_{\rho} \Sigma_{r,s} \to M,
\end{align}
where $\rho := \rho_{r,s}: \Spin_{r,s} \to \GL(\Sigma_{r,s})$ is a complex spinor representation, $\Theta^g:\Spin^g M \to \SO^g M$ is a metric spin structure and $\Spin_{r,s} \subset \GLtp_m$ is the spin group.

\nomenclature[M]{$M$}{a smooth manifold of dimension $m$}
\nomenclature[rs]{$(r,s)$}{signature, $r+s=m$}
\nomenclature[Srs]{$\mathcal{S}_{r,s}(M)$}{space of metrics of signature $(r,s)$}
\nomenclature[Theta]{$\Theta$}{$\GLtp M \to \GLp M$, a topological spin structure}
\nomenclature[SpingM]{$\Theta^g$}{$\Spin^g M \to \SO^g M$, a metric spin structure}
\nomenclature[pigM]{$\pi^g_M$}{$\Sigma^g M \to M$, classical metric spinor bundle }
\nomenclature[tauM]{$\tau_M$}{$\tau_M:TM \to M$, tangent bundle}

\begin{MainThm}[universal spinor bundle, cf. \cref{ThmEmbDirBun}]
	\label{MainThmUSB}
	There exists a finite dimensional vector bundle $\bar \pi^{\Sigma}_{SM}:\bar \Sigma M \to J^1 \pi^{r,s}$ such that for each metric $g \in \mathcal{S}_{r,s}(M)$, there exists $\bar \iota_g$ such that
	\begin{align}
		\label{EqDiagUSBMainResult}
		\begin{split}
			\xymatrix{
				\Sigma^g M
					\ar@{-->}[r]^-{\bar \iota_g}
					\ar[d]^-{\pi^g_M}
				& \bar \Sigma M
					\ar[d]^-{\bar \pi^{\Sigma}_{SM}}
				\\
				M
					\ar[r]^-{j^1(g)}
				& J^1\pi^{r,s}
			}
		\end{split}
	\end{align}
	commutes. Here, $J^1 \pi^{r,s}$ denotes the first jet bundle of $\pi^{r,s}$. Moreover, $\bar \pi^{\Sigma}_{SM}$ carries a connection, a metric and a Clifford multiplication such that $\bar \iota_g$ is a morphism of (generalized) Dirac bundles (see \cref{DefGenDiracBdl}). In addition, $\bar \pi^{\Sigma}_{SM}$ is \emph{natural with respect to spin diffeomorphisms}.
\end{MainThm}

The claims of the last assertions mean that not only the vector bundle structure of any metric spinor bundle $\pi^g_M$ can be recovered from $\bar \pi^{\Sigma}_{SM}$, but also its spinorial connection, its metric and its Clifford multiplication, see \cref{ThmEmbDirBun} for the precise meaning and a proof. We would like to emphasize that the morphism $j^1(g)$ in \cref{EqDiagUSBMainResult} is the $1$-jet of the metric $g$ and that the naturality assertion is formulated with respect to diffeomorphisms and not just isometries, see \cref{SubSectNaturality} for details. 

The rest of this paper is organized as follows: In \cref{SectUSBConstr}, we construct a preliminary version $\pi^{\Sigma}_{SM}:\Sigma M \to \S_{r,s} M$ of $\bar \pi^{\Sigma}_{SM}$ without using the theory of jet spaces and in \cref{SectMetrVertConn}, we show that this almost does the job, see \cref{LemSigmagMFromSigmaMVector}. In \cref{SubSectJetSpaces}, we introduce some basic notions about jet spaces, which will then be applied in \cref{SubSectUsjbConstr} to obtain the bundle $\bar \pi^{\Sigma}_{SM}$ from \cref{MainThmUSB}, see \cref{ThmUSBCExists}. We will then show in \cref{SubSectUnivStructsUsjb} that the additional structures constructed in \cref{SectMetrVertConn} descend to the jet bundle and give a proof of \cref{ThmEmbDirBun}, which directly implies \cref{EqDiagUSBMainResult}. In \cref{SubSectNaturality} we give the claim that the universal spinor bundle is natural a precise category theoretic meaning and a proof, see \cref{ThmUSjBisFunctor}. We also show that there is no way to define a spinor bundle like \cref{EqDefSpinorBundleClassic} that is natural under spin diffeomorphisms, see \cref{ThmNoClassicalUSB}. 

In \cref{SectEinsteinDirac}, we draw the following conclusions from that.

\begin{MainThm}[Universal Dirac operator and Einstein-Dirac-Maxwell equation]$ $ 
	\begin{enumerate}
		\item
			There exists a first-order differential operator on the universal spinor bundle, which is natural under spin diffeomorphisms and corresponds to the usual Dirac operator under compositions with the jet of the metric (and analogously for the Einstein-Dirac-Maxwell-equation, see \cref{ThmDiracNaturalDiffOp} resp. \cref{ThmEDMNaturalDiffOp} for details.)
		\item
			In the compact Riemannian case, the solutions of the Einstein-Dirac-Maxwell equation can be expressed as critical points of a functional whose Lagrangian density factors through the second jet space of a \emph{natural finite dimensional fibre bundle}, see \cref{ThmCritEinstDiracFried} for details.
		\item
			In the Lorentzian case, there exists a maximal Cauchy development for the Einstein-Dirac equation, see \cref{ThmEDMCauchyExists}.
	\end{enumerate}
\end{MainThm}
\vspace{0.5em}

Since the space of sections of a finite dimensional bundle has a nice topology, this makes the equations accessible to techniques like the Palais-Smale condition and Morse theory. As the whole construction is natural in the category of spin manifolds with topological spin structures and spin diffeomorphisms, we obtain a \emph{spin-topological field theory}, i.e. a well-defined action of the group of spin diffeomorphisms on the solution space, whose moduli space could be further examined.

The problem how to deal with the technical issues resulting from the fact that the spinor bundle depends on the metric has been approached in various ways. In the Riemannian case, our construction in \cref{DefusbFun} below agrees with the classical construction in \cite{BourgGaud}; in fact our construction is inspired by this very article. Its main ingredient is the partially defined \emph{vertical connection}, see \cref{DefMetricConnectionUSBnonJet}, which is already enough to construct identification isomorphisms $\beta_{g,h}:\Sigma^g M \to \Sigma^h M$, see \cref{Defbetagh}. In the Riemannian case and in its non-jetted version, the universal spinor bundle $\pi^{\Sigma}_{SM}$ has recently become of particular importance in a study of the \emph{spinorial energy functional}, c.f. \cite{AmmannSpinorial, AmmannSpinorialSurface}. One can also construct the identification isomorphisms $\beta_{g,h}$ without the partial connection, see \cite{MaierGen}, and show that these give the union $L^2(\Sigma M) := \coprod_{g \in \mathcal{S}_{m,0}(M)}{L^2(\Sigma^g M)}$ the structure of a Hilbert bundle over the Riemannian metrics, see \cite[Chapter 4]{NiknoDiss}. But since the base and the fiber of that bundle are both infinite dimensional, the space of sections of this bundle does not have a canonical Fréchet space topology, which makes it difficult to do calculus in this space. This becomes much easier in the space of sections of $\pi^{\Sigma}$ from \cref{EqUniSpinDouble} respectively $\bar \pi^{\Sigma}$ from \cref{EqDiagUSBC}. Another conceptually different approach described in \cite{BaerGaudMor} rests on \emph{generalized cylinders}, which reduce the problem of finding identification isomorphisms between the spinor bundles to the problem of finding a natural path connecting the two metrics, which unfortunately is not always possible (see \cite[Sect. 9]{BaerGaudMor} for a detailed discussion of the Lorentz case). In particular in case of spacetimes one can also use a trivial bundle as a spinor bundle and recover the metric from the operator, see \cite{finsteru22}. However, to study the initial value problem of the Einstein-Dirac-Maxwell equation, we need a stronger notion of naturality of the situation. We will discuss this problem in \cref{SectEinsteinDirac}.

\textbf{Acknowledgements.} The authors would like to thank Bernd Ammann and Felix Finster for interesting discussions and suggestions. 

\vspace{6em}

\section{A universal spinor bundle}

In this section, we review the classical approach by Bourguignon and Gauduchon in \cite{BourgGaud}. We reformulate and slightly generalize their results to metrics of arbitrary signature and show that the resulting vector bundle already carries a lot of additional structure.

\subsection{Construction of the bundle}
\label{SectUSBConstr}

It turns out that the construction of the universal spinor bundle is easier, if we first consider the case of a real vector space. So, let $V$ be an oriented real $m$-dimensional vector space, $\S_{r,s} V$  be the space of inner products of signature $(r,s)$ and $\GLp V$ be the oriented bases on $V$. We denote the action $a$ of a matrix $A \in \GLp_m$ on a basis $b \in \GLp V$ by $ a(b,A) := b.A := b'$, where $b_j' = A^i_j b_i$. We recall that for any metric $g \in \S_{r,s} V$, there exists a \emph{positive pseudo-orthonormal basis $b$}, i.e. a $b \in \GLp V$ satisfying $g(b_i, b_j) = \delta_{ij} \varepsilon_i$, $1 \leq i,j \leq m$, where $\varepsilon_i = +1$ for $1 \leq i \leq r$ and $\varepsilon_i = -1$ for $r+1 \leq i \leq m$. Conversely, any positive basis $b$ of $V$ determines a metric $g_b$ by declaring $b$ to be a pseudo-orthonormal basis, i.e. by setting $g_b(b_i,b_j) := \delta_{ij} \varepsilon_i$. If $b' = b.A$ for some matrix $A \in \GLp_m$, then $g_b = g_{b'}$ if and only if $A \in \SO_{r,s}$. In the language of bundle theory, this fact can be rephrased as follows.

\begin{Lem}
	\label{LemKappaVect}
	\nomenclature[kappa]{$\kappa^V$}{$\GLp V  \to \S_{r,s} V$}
	There is a commutative diagram
	\begin{align}
		\label{EqDefKappaPiPhi}
		\begin{split}
			\xymatrix{
				\GLp V
					\ar[r]^-{\kappa^V}
					\ar[d]^-{q}
				& \S_{r,s}V 
				\\
				\GLp V / \SO_{r,s},
					\ar@{-->}[ur]^-{\varphi}
			}
		\end{split}
	\end{align}
	where, $\kappa^V(b):=g_b$ and $q$ is the canonical quotient map. Moreover, $q$ and $\kappa^V$ are smooth principal $\SO_{r,s}$-bundles and $\varphi$ is a diffeomorphism. In the Riemannian case, this bundle is globally trivial.
\end{Lem}

\begin{Rem}
	\label{DefKappaTildeV}
	Let $\vartheta_{r,s}: \Spin_{r,s} \to \SO_{r,s}$ be the non-trivial $2:1$-cover and $\theta: \GLtp V \to \GLp V$ be the universal cover of $\GLp V$. Then we can extend \cref{EqDefKappaPiPhi} to 
	\begin{align*}
		\xymatrix{
			\GLtp V
				\ar[r]^-{\theta}_{2:1}
				\ar[d]
			& \GLp V
				\ar[r]^-{\kappa^V}
				\ar@{->>}[d]
			& \S_{r,s} V
			\\
			\GLtp V / \Spin_{r,s}
				\ar[r]^-{\cong}
			& \GLp V / \SO_{r,s}.
				\ar@{-->}[ur]_-{\cong }
		}
	\end{align*}
	We set $\tilde \kappa^{V} := \kappa^{V} \circ \theta$.
\end{Rem}

\nomenclature[thetavar]{$\vartheta_{r,s}$}{$2:1$-cover $\Spin_{r,s} \to \SO_{r,s}$}
\nomenclature[thetavar]{$\theta$}{$2:1$-cover $\GLtp_m \to \GLp_m$}
\nomenclature[kappatilde]{$\tilde \kappa^V$}{$\GLtp V \to \S_{r,s} V$}
\nomenclature[rhors]{$\rho_{r,s}$}{fixed choice of spinor representations $\Spin_{r,s} \to \GL(\Sigma_{r,s})$}

Diagram \cref{EqDefKappaPiPhi} easily extends to spin manifolds as follows: Let $\{\rho_{r,s}: \Spin_{r,s} \to \GL(\Sigma_{r,s})\}_{r,s \in \N}$ be a fixed choice of $\Spin_{r,s}$-representations. The construction of $\kappa^V$ and $\tilde \kappa^V$ induces bundle maps
\begin{align}
	\label{EqDefKappaTildeE}
	\kappa^M: \GLp M \to \S_{r,s} M, &&
	\tilde \kappa^M: \GLtp M \to \S_{r,s} M, 
\end{align}
by setting $\kappa^M|_{\GLp_x M} := \kappa^{T_xM}$ and $\tilde \kappa^M|_{\GLtp_x M} := \tilde \kappa^{T_x M}$ for any $x \in M$.

\begin{Def}[universal spinor bundle]
	\label{DefusbFun}
	The map
		\DefMap{\pi^{\Sigma}_{S M}: \Sigma M}{\S_{r,s} M}{{[\tilde b, \sigma]}}{\tilde \kappa^M(\tilde b),}
	where $\Sigma M := \GLtp M \times_{\rho_{r,s}} \Sigma_{r,s}$ is called \emf{universal spinor (vector) bundle}. We obtain the following diagram
	\begin{align}
		\label{EqUniSpinDouble}
		\xymatrixcolsep{4em}
		\begin{split}
			\xymatrix{
				\Sigma M
					\ar[r]^-{\pi^{\Sigma}_{SM}}
					\ar@/_2pc/[rr]^-{\pi^{\Sigma}_{M}}
				& \S_{r,s} M
					\ar[r]^-{\pi^{r,s}}
				& M,
			}
		\end{split}
	\end{align}
	where $\pi^{\Sigma}_M := \pi^{r,s} \circ \pi^{\Sigma}_{SM}$. The map $\pi^{\Sigma}_{M}$ is called \emf{universal spinor (fiber) bundle} of $M$. Its sections are called \emf{universal spinor fields}. 	
\end{Def}

\nomenclature[piSigmaSM]{$\pi^{\Sigma}_{SM}$}{$\Sigma M \to \S_{r,s}M$, universal spinor (vector) bundle}
\nomenclature[piSigmaM]{$\pi^{\Sigma}_M$}{$\Sigma M \to M$, universal spinor (fiber) bundle}

\begin{Rem}[universal spinor fields and spinor fields]
	\label{LemUniversalSpinorFields} 
	For any universal spinor field $\Phi \in \Gamma(\pi^{\Sigma}_M)$, setting $g := g_{\Phi} := \pi^{\Sigma}_{SM} \circ \Phi \in \mathcal{S}_{r,s}(M)$ and $\varphi := \varphi_{\Phi} := \Phi \in \Gamma(\pi^g_M)$ yields a metric on $M$ and an associated spinor field in $\Gamma(\pi^g_M)$. Conversely, we can identify every spinor field $\varphi \in \Gamma(\pi^g_M)$ with a universal spinor field $\Phi := (g, \varphi)$.
\end{Rem}

\begin{Rem}[notation for space of sections]
	For any fiber bundle $\pi_P: P \to M$, we denote its space of sections by $\Gamma(\pi_P)$ instead of $\Gamma(P)$, because we will deal frequently with doubly fibered bundles as in \cref{EqUniSpinDouble} where the projection will be important. We denote by $\tau_P:TP \to P$ the tangent space of the bundle and by $\tau^v_P: T^v P \to P$ its vertical subbundle.
\end{Rem}

\nomenclature[tauv]{$\tau^v_P$}{vertical tangent bundle of $P$}

\subsection{Existence of universal structures on the bundle}
\label{SectMetrVertConn}
We show that the principal $\SO_{r,s}$-fiber bundle $\kappa^V:\GLp V \to \S_{r,s} V$ carries much more structure than just being a fiber bundle. 

\begin{Rem}[metric]
	\label{RemMetricRn}
	\nomenclature[Symrs]{$\Sym_{r,s}$}{$(r,s)$-symmetric matrices}
	\nomenclature[Asymrs]{$\Asym_{r,s}$}{$(r,s)$-anti-symmetric matrices}
	Let $I_{r,s} \in \R^{m \times m}$ be the diagonal matrix with $r$ entries of $+1$ followed by $s$ entries of $-1$ and set $I_m := I_{m,0}$. For any matrix $A \in \R^{m \times m}$, we set $A^\dagger := I_{r,s}A^TI_{r,s}$, denote by $\langle A, B\rangle := \tr(A^\dagger B)$ the metric on $\R^{m \times m}$ and define $\Sym_{r,s} := \{A \in \R^{m \times m} \mid A^{\dagger} = A \}$ and $\Asym_{r,s} := \{A \in \R^{m \times m} \mid A^{\dagger} = -A \}$. Recall that the natural decomposition 
	\begin{align}
		\label{EqDecRnAsym}
		\R^{m \times m} = \Asym_{r,s} \oplus \Sym_{r,s},
	\end{align}
	is orthogonal with respect to $\langle \_, \_ \rangle$. 
\end{Rem}

In a bundle theoretic language, we even get the following.

\begin{Lem}[natural metric and connection]
	\label{LemNatMetrConn}
	\nomenclature[TvGLpV]{$T^v \GLp V$}{vertical distribution on $\GLp V$}
	\nomenclature[ThGLpV]{$T^h \GLp V$}{horizontal distribution on $\GLp V$}
	For any two bases $b,b' \in \GLp V$, let $\tau_b(b') \in \R^{m \times m}$ be the coordinate matrix defined by $b_j' = \tau_b(b')^i_j b_i$, $1 \leq j \leq m$. For any $X,Y \in T_b \GLp V$, we set $\langle X, Y \rangle_b := \langle d \tau_b X, d \tau_b Y \rangle$. Then the $\langle \_, \_ \rangle_b$ assemble to a pseudo-Riemannian metric on $\GLp V$ such that $\SO_{r,s}$ acts by isometries. In particular, 
	\begin{align}
		\label{EqBGHD}
		\forall b \in \GLp V: T^v_b \GLp V := \ker d_b \kappa^V , && T^h \GLp V := (T^v_b \GLp V)^{\perp},
	\end{align}
	defines an orthogonal decomposition such that $T^h \GLp V$ is a connection on $\GLp V$. This decomposition satisfies 
	\begin{align*}
		d_b \tau_b (T^v_b \GLp V) = \Asym_{r,s},
		&&  d_b \tau_b (T^h_b \GLp V) = \Sym_{r,s},
	\end{align*}
	and thus corresponds to \cref{EqDecRnAsym}.
\end{Lem}

An elementary calculation shows that this construction behaves nicely with respect to morphisms.

\begin{Lem}
	\label{LemNatMetricIsom}
	Let $I:V \to W$ be an orientation-preserving isomorphism of oriented vector spaces. Then
	\begin{align*}
		\xymatrixcolsep{3.5em}
		\xymatrix{
			\GLp V
				\ar[d]^-{\kappa^V}
				\ar@{-->}[r]
			&\GLp W
				\ar[d]^-{\kappa^W}
			\\
			\S_{r,s} V
				\ar[r]^-{(I^{-1})^*}
			&\S_{r,s} W
		}
	\end{align*}
	commutes and the top row is an $\SO_{r,s}$-equivariant isometric diffeomorphism. 
\end{Lem}

\begin{Def}
	\label{DefBtpV}
	Let $\Theta^V: \GLtp V \to \GLp V$ be a universal cover. It follows automatically that $\GLtp V$ is a pseudo-Riemannian manifold such that $\Theta^V$ is a local isometry. Therefore, we obtain a connection on $\GLtp V$ by taking the orthogonal complements again. In this manner, $\tilde \kappa^V:\GLtp V \to \S_{r,s} V$ becomes a principal $\Spin_{r,s}$-bundle with a pseudo-Riemannian metric and a connection
	\begin{align}
		\label{EqBGHDTilde}
		T \GLtp V = T^v \GLtp V \oplus T^h \GLtp V.
	\end{align}
	The connection $T^h \GLtp V$ is called the \emf{Bourguignon-Gauduchon horizontal distribution}.
\end{Def}

\nomenclature[TvGLtpV]{$T^v \GLtp V$}{vertical distribution on $\GLtp V$}
\nomenclature[ThGLtpV]{$T^h \GLtp V$}{horizontal distribution on $\GLtp V$}

These constructions carried out on a vector space can be generalized to a spin manifold $(M, \Theta)$ as follows.

\begin{Def}[vertical connection]
	\label{DefMetricConnectionUSBnonJet}
	For any $x \in M$, we define
	\begin{align}
		\label{EqDefVertDist}
		T^{vv}_x \GLp M := T^v \GLp (T_x M), &&
		T^{vh}_x \GLp M := T^h \GLp (T_x M), 
	\end{align}
	and analogously for $\GLtp M$. The resulting decomposition $T^v \GLp M = T^{vv} \GLp M \oplus T^{vh} \GLp M$	is called a \emf{vertical distribution} on $\kappa^M: \GLp M \to \S_{r,s} M$ (and analogously on $\tilde \kappa^M:\GLtp M \to \S_{r,s} M$).  We denote by
	\begin{align}
		\label{EqDefTildeNablaE}
		\bm{\nabla}:\Gamma(\tau^v_{\S_{r,s}M}) \times \Gamma(\pi^{\Sigma}_{SM}) \to \Gamma(\pi^{\Sigma}_{SM})
	\end{align}
	the induced \emf{vertical connection}, i.e. the connection induced on $\pi^{\Sigma}_{SM}$ from \cref{EqUniSpinDouble} that is only defined for all directions in the vertical space $\tau^v_{\S_{r,s}M}:T^v \S_{r,s}M \to M$.
\end{Def}

\nomenclature[nablabold]{$\bm{\nabla}$}{vertical spinorial connection}

In the Riemannian case, this connection agrees with the original construction in \cite{BourgGaud}. We now show that this bundle carries a lot more structure than just being a bundle with a partially defined connection. We denote by $\tau^{r,s}_M := (\pi^{r,s})^*(TM) \to \S_{r,s}M$ the pullback of the tangent bundle of $M$ to $\S_{r,s} M$ by $\pi^{r,s}$.

\nomenclature[taursM]{$\tau^{r,s}_M$}{$(\pi^{r,s})^*(TM) \to \S_{r,s}M$}
\nomenclature[eta]{$\bm{\eta}$}{universal metric on $\pi^{\Sigma}_{SM}$}
\nomenclature[m]{$\mathfrak{m}$}{universal Clifford multiplication on $\pi^{\Sigma}_{SM}$}

\begin{Def}[universal metric and Clifford multiplication]
	\label{DefUniMetricCliffMult}
	We define
	\begin{align}
		\label{EqDefEta}
		\bm{\eta}(\phi, \phi') & :=  g_x(\phi, \phi'), \\
		\label{DefCliffordUniv}
		\mathfrak{m}(X^* \otimes \phi) & := X^* \bullet \phi := V \cdot_{g_x} \phi,
	\end{align}
	where $\phi, \phi'  \in \Sigma M|_{g_x}$, $X^* = (g_x, V) \in (\pi^{r,s})^*(TM)$, $V \in T_x M$. The map $\bm{\eta}$ is called \emf{universal spinorial metric} and $\mathfrak{m}$ is called \emf{universal Clifford multiplication} on $\pi^{\Sigma}_{SM}:\Sigma M \to \S_{r,s}M$.
\end{Def}

One might wonder to what extent $\bm{\nabla}$, $\bm{\eta}$ and $\mathfrak{m}$ are compatible with each other. To see that, first recall how the analogous structures on the metric spinor bundle $\pi^g_M: \Sigma^g M \to M$ are related to each other.

\begin{Rem}[Dirac structure on the metric spinor bundle]
	\label{RemDiracStructSpingM}
	Recall that a metric spinor bundle $\pi^g_M: \Sigma^g M \to M$ carries a canonical structure as a \emf{Dirac bundle}, i.e. there exists an extension of the Levi-Civita connection to the spinor bundle, also denoted by $\nabla^g: \Gamma(\pi^g_M) \to \Gamma(\tau^*_M \otimes \pi^g_M)$, an extension of the metric $g$ to the spinor bundle, also denoted by $g$, and a \emf{Clifford multiplication}, i.e. a morphism of real vector bundles $\mathfrak{m}^g: TM \otimes_{\R} \Sigma^g M \to \Sigma^g M$, $V \otimes \psi \mapsto V \cdot_g \psi$. In addition, these structures satisfy the following compatibility relations (see for instance \cite[Th. 1.12]{Baum1981}):
	\begin{align}
		\label{EqNablaClifford}
		- 2 g(V,W) \psi  &= V \cdot_g W \cdot_g \psi  + W \cdot_g V \cdot_g \psi, \\
		\label{EqNablagMetric}
		\nabla^g_V g(\psi, \psi') &= g(\nabla^g_V \psi, \psi') + g(\psi, \nabla^g_V \psi'), \\ 
		\label{EqNablagCliffordComp}
		\nabla^g_V(W \cdot_g \psi) &= \nabla^g_V W \cdot_g \psi + W \cdot_g \nabla^g_V \psi, \\
		\label{EqCliffordHermitian}
		g(V \cdot_g \psi, \psi') & = (-1)^{s+1} g(\psi, V \cdot_g \psi'), 
	\end{align}
	where $\psi, \psi' \in \Gamma(\pi^g_M)$, $V, W \in \Gamma(\tau_M)$. 
\end{Rem}

\nomenclature[nablag]{$\nabla^g$}{(spinorial) Levi-Civita connection}
\nomenclature[mg]{$\mathfrak{m}^g$}{Clifford multiplication w.r.t. $g$}

The universal structures $\bm{\nabla}$, $\bm{\eta}$ and $\mathfrak{m}$ satisfy compatibility relations similar to \cref{EqNablagMetric,EqNablagCliffordComp,EqCliffordHermitian}. In order to be able to precisely formulate them, we need the following notion.

\begin{Def}[universal metric and vertical Levi-Civita connection]
	The metric $\mathbf{g}$ defined by 
	\begin{align}
		\label{EqDefUnivMetricVertical}
		\mathbf{g}(X^*,X^*) := g_x(X,X),
	\end{align}
	where $X^* = (g_x, V) \in (\pi^{r,s})^*(TM)$ is called \emf{universal metric}. Let $\pi^+:\GLp M \to M$ the canonical projection. Via the pullback to $\bar \kappa^M:(\pi^{r,s})^*(\GLp M) \to \S_{r,s}M$ along $\pi^{r,s}:\S_{r,s}M \to M$ we also obtain a vertical connection on $\bar \kappa^M$ from the vertical connection on $\kappa^M$ from \cref{DefMetricConnectionUSBnonJet}. Let $\gamma:\GLp_m \to \GL(\R^m)$ be the standard representation (given by matrix multiplication). Recall that $\GLp M \times_{\gamma} \R^m = TM$, so $(\pi^{r,s})^*(TM) = (\pi^{r,s})^*(\GLp M) \times_{\gamma} \R^m$ and therefore, we obtain a vertical connection on $\tau^{r,s}_M:(\pi^{r,s})^*(TM) \to \S_{r,s} M$, which is denoted by
	\begin{align}
		\label{EqDefNablaM}
		\nabla: \Gamma(\tau^v_{\pi^{r,s}}) \times \Gamma(\tau^{r,s}_M) \to \Gamma(\tau^{r,s}_M).
	\end{align}
	We call $\nabla$ the \emf{vertical universal Levi-Civita connection}.
\end{Def}

\nomenclature[gbold]{$\bm{g}$}{universal metric on $(\pi^{r,s})^*(TM)$}
\nomenclature[nabla]{$\nabla$}{vertical universal Levi-Civita connection}

\begin{Lem}[properties of universal structures]
	\label{LemPropUnivStruct}
	The universal structures satisfy the compatibility conditions 
	\begin{align}
		\label{EqBullletIsClifford}
		-2 \bm{g}(X^*, X^*) &=  X^* \bullet Y^* \bullet \psi + Y^* \bullet X^* \bullet \psi \\
		\label{EqTildeNablaMetric}
		\bm{\nabla}_X( \bm{\eta}(\phi, \phi')) & = \bm{\eta}(\bm{\nabla}_X \phi, \phi') + \bm{\eta}(\phi, \bm{\nabla}_X \phi'), \\
		\label{EqTildeNablaClifford}
		\bm{\nabla}_{X} (Y^* \bullet \phi) & = \nabla_X Y^* \bullet \phi + Y^* \bullet \bm{\nabla}_{X} \phi, \\
		\label{EqEtaCliffordHermitian}
		\bm{\eta}(X^* \bullet \phi, \phi') & = (-1)^{s+1} \bm{\eta}(\phi, X^* \bullet \phi).
	\end{align}
	where $X^*, Y^* \in \Gamma(\tau^{r,s}_M)$, $X \in \Gamma(\tau^v_{\S_{r,s} M})$, $\phi, \phi' \in \Gamma(\pi^{\Sigma}_{SM})$. 
\end{Lem}

\begin{Prf}
	\cref{EqBullletIsClifford} follows directly from \cref{EqNablaClifford} and \cref{EqDefUnivMetricVertical}.
	
	By construction \cref{EqDefEta}, the metric $\bm{\eta}$ agrees pointwise with the spinorial metric, which is defined via the invariant metric on the spinor space $\Sigma_{r,s}$ in the spinor representation $\rho:\Spin_{r,s} \to \GL(\Sigma_{r,s})$. The representation $\rho$ acts on $\Sigma_{r,s}$ by isometries. Since the vertical connection on $\pi^{\Sigma}_{SM}$ is induced by the vertical connection on $\tilde \kappa^M$, the resulting vertical connection is metric. For connections this can be found for instance in \cite[Satz 3.13]{Baum} and the same proof holds for vertical connections. This proves \cref{EqTildeNablaMetric}.
	
	To see \cref{EqTildeNablaClifford}, we just notice that $\nabla \otimes \bm{\nabla}$ is a connection on $(\pi^{r,s})^*(TM) \otimes \Sigma M$ and $\mathfrak{m}$ is parallel with respect to this connection. This results from the fact that the ordinary Clifford multiplications $\mathfrak{m}^g$ are defined via the representation $\rho$. 
	
	Finally, \cref{EqEtaCliffordHermitian} follows from the definition of $\bm{\eta}$, $\bullet$ and \cref{EqCliffordHermitian}. 
\end{Prf}

The relationship between the usual Dirac structure on the metric spinor bundle $\pi^g_M:\Sigma^g M \to M$ from \cref{RemDiracStructSpingM} and the universal structures on $\pi^{\Sigma}_{SM}:\Sigma M \to \S_{r,s} M$ is as follows. 

\begin{Lem} 
	\label{LemSigmagMFromSigmaMVector}
	For any metric $g$, there exists a morphism $I_g$ of vector bundles such that
	\begin{align}
		\label{EqUnivSpinor}
		\begin{split}
			\xymatrix{
				\Sigma^g M
					\ar@/^2pc/[drr]^-{\iota_g}
					\ar@/_2pc/[ddr]_-{\pi^g_M}
					\ar@{-->}[dr]^{I_g}
				\\
				& g^*\Sigma M
					\ar[r]^-{g^*}
					\ar[d]
				& \Sigma M
					\ar[d]^{\pi^{\Sigma}_{SM}}
				\\
				&M
					\ar[r]^-{g}
				&\S_{r,s} M.
			}
		\end{split}
	\end{align}
	commutes. In addition, $I_g$ is an isometric isomorphism with respect to the spinorial metric on $\pi^g_M$ and $g^* \bm{\eta}$ and it is compatible with the Clifford multiplications $\mathfrak{m}^g$ and $g^* \mathfrak{m}$.
\end{Lem}

\begin{Prf}
	For any $x \in M$, $\psi \in \Sigma^g_x M$, we have $(\pi^{\Sigma}_{SM} \circ \iota_g)(\psi)=g_x=(g \circ \pi^g_M)(\psi)$, thus the desired isomorphism exists by the universal property of the pullback. To see the second claim, we check that
	\begin{align*}
		(g^* \bm{\eta})(I_g(\psi), I_g(\psi')) & = \bm{\eta}(\iota_g(\psi), \iota_g(\psi')) = g(\psi, \psi'), \\
		(g^* \mathfrak{m})(V \otimes I_g(\psi)) &= \mathfrak{m}((g,V) \otimes \iota_g(\psi)) = V \cdot_g \psi = \mathfrak{m}^g(V \otimes \psi), 
	\end{align*}
	where $\psi, \psi' \in \Gamma(\pi^g_M)$, $V \in TM$.
\end{Prf}

Notice that $I_g$ is not compatible with the connections, since we would have to check that $(g^* \bm{\nabla})_V(I_g(\psi)) = \bm{\nabla}_{dg V}(\iota_g(\psi)) = \nabla^g_{dgV} \psi$, which makes no sense, since $dg V$ is not vertical, so $\bm{\nabla}_{dgV}$ is not defined. We will discuss this issue in the following.

\subsection{Jet spaces}
\label{SubSectJetSpaces}

The fact that the connections $\bm{\nabla}$ and $\nabla$ are defined only on the vertical part of $\pi^{r,s}:\S_{r,s}M \to M$ is a bit unpleasant, especially since this prevents us from including them into the compatibility assertion of \cref{LemSigmagMFromSigmaMVector}. In this and the next subsection, we show how to obtain a full connection by passing to the first jet space. We consider the first jet bundle, since the spinorial Levi-Civita connection on a spinor bundle depends only on the $1$-jet of the  metric. We recall the definition of a jet space, see \cite{saunders} for a more comprehensive introduction to the topic.

\begin{Def}[jet bundle]
	Let $\pi_P: P \to M$ be any fiber bundle. For any $x \in M$, denote by $\Gamma_x(\pi_P)$ the space of sections defined on a local neighborhood near $x$. Two such sections $s_1, s_2$ \emf{have the same $1$-jet at $x$}, if $s_1(x) = s_2(x) \in P$ and $d s_1 |_{T_xM} = d s_2 |_{T_xM}$. The equivalence class $j^1_x(s)$ of a local section $s \in \Gamma_x(\pi_P)$ is the \emf{$1$-jet of $s$ at $x$}. The set $J^1 \pi_P := \{j^1_x(s) \mid x \in M, s \in \Gamma_x(\pi_P) \}$ is the \emf{first jet space of $\pi_P$}. The space $J^1 \pi_P$ comes along with the two canonical projections $\pi_0 := j^1_0 \pi_P:J^1 \pi_P \to M$, $j^1_x s \mapsto x$, and $\pi_{1,0} := j^1_{1,0}\pi_P:J^1 \pi_P \to P$, $\quad j^1_x s \mapsto  s(x)$, called the \emf{source} respectively \emf{target projection}.
\end{Def}

\nomenclature[J1P]{$J^1 \pi_P$}{first jet bundle of $\pi_P$}

It is well known that $J^1 \pi_P$ is a smooth manifold and $\pi_{1,0}$ is even an affine bundle. One can also define higher jet bundles $J^k \pi_P$ consisting of equivalence classes of sections that agree up to the $k$-th derivative.

\begin{Def}
	Let $\pi_P: P \to M$ be a bundle, $x \in M$, $s \in \Gamma_x(\pi_P)$ and $V \in T_pM$. Then $(j^1_x(s), d_x s (V)) \in \pi_{1,0}^*(TP)$ is called the \emf{holonomic lift of $V$} at $s$. 
\end{Def}

The following theorem asserts that one can decompose $\pi_{1,0}^* (TP)$ into vertical vectors and holonomic lifts.

\begin{Thm}[\protect{\cite[Thm. 4.3.2]{saunders}}]
	\label{ThmUSBJet}
	Let $\pi_{P}:P \to M$ be a fiber bundle and $\pi_{1,0}:J^1 \pi_P \to P$ be the target projection. Recall that $\tau^v_P: T^v P \to P$ denotes the vertical tangent bundle of $P$. At any point $j^1_x(s)$, $s \in \Gamma(\pi_P)$, there exists a natural decomposition
	\begin{align}
		\label{EqJetBundleDecomposition}
		\pi_{1,0}^*(TP)|_{j^1_x(s)} = \pi_{1,0}^*(T^v P)|_{j^1_x(s)} \oplus \underbrace{(j^1_x(s),ds(T_x M))}_{\in \pi_{1,0}^*(TP)|_{j^1_x(s)}}.
	\end{align}
	This decomposition is well-defined, i.e. it does not depend on the choice of $s$ for a given $j^1_x(s)$. 
\end{Thm}

\begin{Rem}  
	\label{RemTangentialJet}
	Consider the canonical map $\psi:TJ^1\pi_P \to \pi_{1,0}^*(TP)$, $ \bar X \mapsto (\tau_{J^1 P}\bar X, d \pi_{1,0} \bar X)$. For any $\bar X \in T_{j^1_x(s)} J^1 \pi_P$, we obtain a decomposition of $\psi(\bar X)$ as in \cref{EqJetBundleDecomposition}. Setting $X := d\pi_{1,0} \bar X$, this is can be written explicitly as 
	\begin{align}
		\label{EqSaundersDec}
		X = X^v \oplus X^h \in T^v P \oplus d_xg(T_xM),
	\end{align}
	i.e. there exists $X_h \in T_xM $ such that $X^h = d_xg(X_h)$. Notice that $X_h = d_{j^1_x(s)} \pi_0 \bar X$.
\end{Rem}

\subsection{Construction of the jet version}
\label{SubSectUsjbConstr}
Before we can state \cref{ThmUSBCExists}, we recall some facts about pullbacks.

\begin{Rem}[description of pull-back connections]
	\label{RemPullback}
	Let $\pi_E: E \to M$ be a vector bundle and $f:N \to M$ be smooth. Then we denote by $f^* E \to N$ the pullback, i.e.
	\begin{align*}
		\xymatrix{
			f^*E
				\ar[r]^-{f^*}
				\ar[d]^{\bar \pi_E}
			&E
				\ar[d]^{\pi_E}
			\\
			N
				\ar[r]^-{f}
			&M.
		}
	\end{align*}
	Recall that $f^*E = \{ (y,e) \in N \times E \mid f(y)=\pi_E(e)\}$. Therefore, a section $s \in \Gamma(\pi_E)$ determines a section $(\id_N, s \circ f) \in \Gamma(\bar \pi_E)$. Conversely, a section $\bar s = (\bar s_1, \bar s_2) \in \Gamma(\bar \pi_E)$ always satisfies $\id_N = \bar \pi_E \circ \bar s = \bar s_1$ and $(\pi_E \circ f^*)(\bar s) = \pi_E \circ  \bar s_2 = f$ and is therefore already determined by $s_2$. We call a map $s:N \to E$ such that $\pi_E \circ s = f$ a \emf{section along $f$} and identify the space $\Gamma(\pi_E) \circ f$ of sections along $f$ with with $\Gamma(\bar \pi_E)$. 
	
	In addition, assume that $E$ carries a connection $\nabla$. Then we can define $ \bar D_X(s \circ f) := \nabla_{f_* X} s$ for any $X \in TN$ and $s \in \Gamma(\pi_E)$. The result is an $\R$-bilinear map $\bar D: \Gamma(\tau_N) \times (\Gamma(\pi_E) \circ f) \to (\Gamma(\pi_E) \circ f)$ such that for all $\bar X \in TN$, $s \in \Gamma(\pi_E)$, we have the usual Leibniz rule and $\mathcal{C}^{\infty}$-linearity:
	\begin{align*}
		\forall \beta \in \mathcal{C}^{\infty}(M): \bar D_{\bar X}( \beta s \circ f) 
		&= \bar X(\beta \circ f) (s \circ f) +  (\beta \circ f) \bar D_{\bar X}(s \circ f), \\
		\forall \alpha \in \mathcal{C}^{\infty}(N): \bar D_{\alpha \bar X}(s \circ f) &= \alpha \bar D_{\bar X}(s \circ f).
	\end{align*}
	It is a standard result in differential geometry that any such map $\bar D$ defines a unique connection $\bar \nabla:\Gamma(\tau_N) \times \Gamma(\bar \pi_E) \to \Gamma(\bar \pi_E)$ satisfying
	\begin{align}
		\label{EqConnConstrProp}
		\nabla _{\bar X}(\id_N,s \circ f) = (\id_N, \bar D_{\bar X}(s \circ f)).
	\end{align}
	We sometimes write $s \circ f$ instead of $(\id_N, s\circ f)$, since it makes no distinction if we use the pullback connection on $\Gamma(\bar \pi_E)$.
\end{Rem}

With this in mind, we can continue \cref{EqUniSpinDouble} (drawn vertically) as follows.

\begin{Thm}[universal spinor jet bundle]
	\label{ThmUSBCExists}
	The universal spinor bundle from \cref{EqUniSpinDouble} can be extended to a commutative diagram
	\begin{align}
		\label{EqDiagUSBC}
		\begin{split}
			\xymatrixcolsep{3.5em}
				\xymatrix{
					\bar \Sigma M
						\ar[r]^{F^{\Sigma}}
						\ar[d]^{\bar \pi^{\Sigma}_{SM}}
						\ar@/_3em/[dd]_-{\bar \pi^{\Sigma}_{M}}
					& \Sigma M
						\ar[d]^-{\pi^{\Sigma}_{SM}}
						\ar@/^3em/[dd]^-{\pi^{\Sigma}_{M}}
					\\
					J^1 \pi^{r,s}
						\ar[r]^-{\pi^{r,s}_{1,0}}
						\ar[d]^{\pi_0^{r,s}}
					& \S_{r,s} M
						\ar[d]^{\pi^{r,s}}
					\\
					M
						\ar[r]^-{\id_M}
					& M.
				}
		\end{split}
	\end{align}
	Moreover, the vector bundle $\bar \pi^{\Sigma}_{SM}$ carries a connection $\bm{\bar \nabla}$ satisfying 
	\begin{align}
		\label{EqDefUnivConn}
		F^{\Sigma}(\bm{\bar{\nabla}}_{\bar X} \bar \phi|_{j^1_x(g)}) 
		=  \bm{\nabla}_{X^v} \phi |_{g(x)}  + \nabla^g_{X_h}(\phi \circ g)|_x,
	\end{align}
	where $\bar \phi := (\id, \pi^{r,s}_{1,0} \circ \phi) \in \Gamma(\bar \pi^{\Sigma}_{SM})$, $\phi \in \Gamma(\pi^{\Sigma}_{SM})$, $\bar X \in T_{j^1_x(g)}(J^1 \pi^{r,s})$, $X := d \pi^{r,s}_{1,0} \bar X$ and $X^v$ and $X_h$ are a decomposition as in \cref{EqSaundersDec}. Here, $\pi^{r,s}_{1,0} = j^1_{1,0} \pi^{r,s}$ and $\pi^{r,s}_0 := j^1_0 \pi^{r,s}$.
\end{Thm}

\nomenclature[nablaboldbar]{$\bm{\bar \nabla}$}{connection on universal spinor jet bundle}
\nomenclature[pibarSigmaSM]{$\bar \pi^{\Sigma}_{SM}$}{$\bar \Sigma M \to J^1\pi^{r,s}$, jetted universal spinor (vector) bundle}
\nomenclature[pibarSigmaM]{$\bar \pi^{\Sigma}_{M}$}{$\bar \Sigma M \to M$, jetted universal spinor (fiber) bundle}

\begin{Prf}
	The strategy is to define the connection on $\bar \Sigma M$ using \cref{RemPullback}. We define $\bar \Sigma M := (\pi^{r,s}_{1,0})^* \Sigma M $ and obtain commutativity of \cref{EqDiagUSBC}. Using the notation from the assertion and \cref{EqSaundersDec}, we obtain a decomposition
		\begin{align}
			\label{EqJetDecompositionUSBC}
			X = X^v \oplus X^h \in T^v \S_{r,s} M \oplus d_xg (T_x M),
		\end{align}
	i.e. there exists $X_h \in T_x M$ such that $X^h = d_x g(X_h)$. Recall from \cref{LemUniversalSpinorFields} that $\phi \circ g \in \Gamma(\pi^g_M)$ and that the connection $\nabla^g$ on the spinor bundle at a point $ x \in M$ depends only on $j^1(g)$ at $x$, so the right hand side of \cref{EqDefUnivConn} is well-defined. Therefore, we simply define $\bar D_X \bar \phi$ by the right hand side of \cref{EqDefUnivConn}. To show that this gives a connection on $\bar \kappa^{\Sigma}_{SM}$, it remains only to verify the properties of $\bar D$ enlisted in \cref{RemPullback}.
	It is clear that $\bar D$ is $\R$-bilinear in both arguments. To see the Leibniz rule, let $\beta \in \mathcal{C}^{\infty}(\S_{r,s}M)$ and calculate
	\begin{align*}
		\bar D_{\bar X} ((\beta \phi) \circ \pi^{r,s}_{1,0}) |_{j^1_x(g)}
		&= \bm{\nabla}_{X^v}(\beta \phi)|_{g(x)} + \nabla^g_{X_h}(\beta \phi \circ g) |_x \\
		&= X^v(\beta) \phi|_{g(x)} 
		+ \beta \bm{\nabla}_{X^v}(\phi)|_{g(x)} 
		+ X_h(\beta \circ g) \phi|_x 
		+ (\beta \circ g) \nabla^g_{X_h}(\phi \circ g)|_x\\
		&= ( X^v(\beta) 
		+ dg(X_h)(\beta))|_{g(x)}\phi|_{g(x)} 
		+ \beta \bm{\nabla}_{X^v}(\phi)|_{g(x)}
		+(\beta \circ g) \nabla^g_{X_h}(\phi \circ g)|_x\\
		&= X(\beta)|_{g(x)} \phi|_{g(x)} 
		+ (\beta \circ g)|_x \bar D_{\bar X}(\phi \circ \pi^{r,s}_{1,0})|_{j^1_x(g)} \\
		&= \bar X(\beta \circ \pi^{r,s})|_{j^1_x(g)} \phi|_{g(x)}
		+ (\beta \circ \pi_{1,0}^{r,s})|_{j^1_x(g)} \bar D_{\bar X}(\phi \circ \pi^{r,s}_{1,0})|_{j^1_x(g)}.
	\end{align*}
	To see the $\mathcal{C}^\infty$-linearity, let $\alpha \in \mathcal{C}^{\infty}(J^1 \pi^{r,s})$. We obtain $d\pi^{r,s}_{1,0}(\alpha \bar X) = \alpha d\pi^{r,s}_{1,0} \bar X$ and clearly $dg(\alpha(j^1_p(g)) X_h) = \alpha(j^1_p(g)) X^h $, thus 
	\begin{align*}
		\bar D_{\alpha(j^1_p(g)) \bar X} (\phi \circ \pi^{r,s}_{1,0}) 
		&= \bm{\nabla}_{\alpha(j^1_p(g)) X^v} \phi + \nabla^g_{\alpha(j^1_p(g)) X_h}(\phi \circ g) \\
		&= \alpha(j^1_p(g)) (D_{\bar X} (\phi \circ \pi^{r,s}_{1,0})),
	\end{align*}
	which concludes the proof.
\end{Prf}

\subsection{Existence of universal structures on the jet version}
\label{SubSectUnivStructsUsjb}

We will now show that the structures $\bm{g}$, $\bm{\eta}$, $\mathfrak{m}$ can also be pulled back to the jet bundle $\bar \pi^{\Sigma}_{SM}$ and satisfy compatibility conditions with the connection $\bm{\bar \nabla}$ similar to \cref{LemPropUnivStruct}. This also allows us to formulate a version of \cref{LemSigmagMFromSigmaMVector}, which includes the connection. 

\begin{Def}[universal Dirac structure]
	\label{DefUnivDiracStruct}
	Consider the vector bundle $\bar \pi^{\Sigma}_{SM}: \bar \Sigma M \to J^1 \pi^{r,s}$ from \cref{EqDiagUSBC}. We define
	\begin{align}
		\label{EqDefUniMetricgbar}
		\overline{\mathbf{g}}(\bar X^*, \bar Y^*) &:= g_x(V, W), \\
		\label{EqDefUnivSpinBldMetric}
		\bm{\bar \eta}_{j^1_x(g)}(\bar \phi, \bar \phi') & := \bm{\eta}_{g_x} (\phi, \phi'), \\
		\label{DefCliffordBar}
		\overline{\mathfrak{m}}(\bar X^* \otimes \bar \phi) & := \bar X^* \bbullet \bar \phi := V \cdot_{g_x} \phi,
	\end{align}
	where $j^1_x(g) \in J^1\pi^{r,s}$,  $\bar X^* = (j^1_x(g), V), \bar Y^* = (j^1_x(g), W) \in (\pi^{r,s}_0)^* TM)$, $\bar \phi, \bar \phi' \in \Gamma(\bar \pi^{\Sigma}_{SM})$, $\phi := F^{\Sigma}(\bar \phi)$, $\phi' := F^{\Sigma}(\bar \phi')$. 
\end{Def} 

\nomenclature[gbar]{$\bm{\bar g}$}{universal metric on $(\pi^{r,s}_0)^*(TM)$}
\nomenclature[etabar]{$\bm{\bar \eta}$}{universal metric on $\bar \pi^{\Sigma}_{SM}$}
\nomenclature[mbar]{$\overline{\mathfrak{m}}$}{Clifford multiplication for $\bar \pi^{\Sigma}_{SM}$}
\nomenclature[nablabar]{$\bar \nabla$}{universal Levi-Civita connection}

Of course, one can also pullback the vertical connection $\nabla$ from \cref{EqDefNablaM} to a vertical connection $(\pi^{r,s}_{1,0})^* \nabla$. This vertical connection can be completed to a full connection in complete analogy to \cref{EqDefUnivConn}.

\begin{Def}[universal Levi-Civita connection] 
	For any $\bar X \in T_{j^1_x(g)} J^1\pi^{r,s}$, $\bar Y^* = (j^1_x(g), V) \in (\pi^{r,s}_0)^*(TM)$, $\bar Y =(g_x, V) \in (\pi^{r,s})^*(TM)$, $V \in T_xM$, we set
	\begin{align*}
		\bar \nabla_{\bar X} \bar Y^* |_{j^1_x(g)}
		:= \nabla_{X^v}{Y^*} 
		+ \nabla^g_{X_h} V ,		
	\end{align*} 
	where $d\pi^{r,s}_{1,0}\bar X=: X = X^v \oplus d_xg X_h$ is decomposed as in \cref{EqSaundersDec}. The connection $\bar \nabla$ is called \emf{universal Levi-Civita connection}. 
\end{Def} 

To see that $\bar \nabla$ is well-defined, we notice that the Levi-Civita connection only depends on the $1$-jet of the metric. So one can proceed as in the proof of \cref{ThmUSBCExists}.

In case $r=m$, the metric $\overline{\mathbf{g}}$ is precisely the \emf{universal Riemannian metric} and $\bar \nabla$ is the universal Levi-Civita connection as considered for instance in \cite{perez}. The equality between $\bar \nabla$ and the connection considered in \cite{perez} is evident from \cite[Thm. 5.1]{perez},  where the uniqueness of natural metric connections on the bundle is shown, and the fact that the connection constructed here is both metric and natural under the action of diffeomorphisms of $M$, see \cref{EqUnivStuctsNatural}.

\begin{Lem}[properties of universal Dirac structures]
	The structures from \cref{DefUnivDiracStruct} satisfy the following compatibility relations: 
	\begin{align} 
		\label{EqbarCliffordRel}
		- 2 \overline{\mathbf{g}}(\bar X^*, \bar Y^*) \bar \phi &= \bar X^* \bbullet \bar Y^* \bbullet \bar \phi + \bar Y^* \bbullet \bar X^* \bbullet \bar \phi , \\
		\label{EqbarNablaMetric}
		\bm{\bar \nabla}_{\bar X} \bm{\bar \eta}(\bar \phi, \bar \phi ') & = \bm{\bar \eta}( \bm{\bar \nabla}_{\bar X} \bar \phi, \bar \phi') + \bm{\bar \eta} (\bar \phi, \bm{\bar \nabla}_{\bar X} \bar \phi'), \\
		\label{EqbarNablaClifford}
		\bm{\bar \nabla}_{\bar X}{(\bar Y^* \bbullet \bar \phi)} &= \bar \nabla_{\bar X} \bar Y^* \bullet \bar \phi + \bar Y^* \bbullet \bm{\bar \nabla}_{\bar X} \bar \phi, \\
		\label{EqbarCliffordHermitian}
		\bm{\bar \eta}(\bar X \bbullet \phi, \phi') &= (-1)^{s+1} \bm{\bar \eta}(\phi, \bar X \bbullet \phi'),
	\end{align}
	where $\bar X \in T J^1 \pi^{r,s}$, $\bar X^*, \bar Y^* \in (\pi^{r,s}_0)^*(TM)$, $\bar \phi, \bar \phi' \in \Gamma(\bar \pi^{\Sigma}_{SM})$.
\end{Lem}

\begin{Prf}	
	The equation \cref{EqbarCliffordRel} follows directly from \cref{EqDefUniMetricgbar,DefCliffordBar,EqBullletIsClifford}.

	To see \cref{EqbarNablaMetric}, we decompose 
	\begin{align*}
			d \pi^{r,s}_{1,0} \bar X =: X = X^v \oplus X^h \in T^v \S_{r,s} M \oplus d_p g (T_pM),
			&& X^h = dg(X_h),
			&& X_h \in T_p M, 
	\end{align*}
	as in \cref{EqJetDecompositionUSBC}. We also assume that $\bar \phi = \pi^{r,s}_{1,0} \circ \phi$ for some $\phi \in \Gamma(\pi^{\Sigma}_{SM})$ and analogously for $\bar \phi '$. By definition
	\begin{align*}
		\bm{\bar \nabla}_{\bar X} \bm{\bar \eta}(\bar \phi, \bar \phi')
		& = \bar X( \bm{\bar \eta} (\bar \phi, \bar \phi')) 
		= X(\bm{\eta}(\phi, \phi') \circ \pi^{r,s}_{1,0}) 
		= d \pi^{r,s}_{1,0} \bar X( \bm{\eta} (\phi, \phi')) \\
		& = X^v(\bm{\eta}(\phi, \phi'))+ X^h (\bm{\eta} (\phi, \phi')). 
	\end{align*}
	Now, for the vertical part, \cref{EqTildeNablaMetric} implies
	\begin{align*}
		X^v(\bm{\eta}(\phi,\phi'))
		= \bm{\eta}(\bm{\nabla}_{X^v} \phi, \phi') + \bm{\eta}(\phi, \bm{\nabla}_{X^v}\phi').
	\end{align*}
	Since $\phi \circ g \in \Gamma(\pi^g_M)$, we obtain for the horizontal part
	\begin{align*}
		X^h(\bm{\eta}(\phi,\phi')
		& = dg(X_h)(\bm{\eta}(\phi, \phi'))
		= X_h(\bm{\eta}(\phi, \phi') \circ g)
		=X_h(g(\phi \circ g, \phi' \circ g)) \\
		& =g( \nabla^g_{X^h} (\phi \circ g), \phi' \circ g) 
		+ g(\phi \circ g, \nabla^g_{X_h}(\phi' \circ g))
	\end{align*}
	Since,
	\begin{align*}
		\bm{\eta}(\bm{\nabla}_{X^v} \phi, \phi') + g(\nabla^g_{X_h}(\phi \circ g|_x, \phi' \circ g))
		& = \bm{\bar \eta}((\id, \bm{\nabla}_{X^v} \phi), \phi')
		+ \bm{\eta}(\nabla^g_{X_h}(\phi \circ g)_x, \phi'|_{g_x}) \\
		&=\bm{\bar \eta}(\bm{\bar \nabla}_{\bar X} \bar \phi, \bar \phi'),
	\end{align*}
	this implies \cref{EqbarNablaMetric}.
	
	To see \cref{EqbarNablaClifford}, we consider $Y^* = (j^1_x(g), V) \in (\pi^{r,s})^*(TM)$, $V \in T_xM$, $\bar Y^* = (g_x, V) \in (\pi^{r,s}_0)^*(T_xM)$. Using \cref{EqNablagCliffordComp,EqTildeNablaClifford}, we obtain
	\begin{align*}
		\bm{\bar \nabla}_{\bar X}(\bar Y^* \bbullet \bar \phi)
		& = \bm{\nabla}_{X^v}(Y^* \bullet \phi) \circ g + \nabla^g_{X_h}((V \cdot \phi) \circ g) \\
		& = \nabla_{X^v}Y^* \bullet \phi 
		+ Y^* \bullet \bm{\nabla}_{X^v}{\phi}
		+ (\nabla^g_{X^h}V) \cdot \phi \circ g
		+ V \cdot \nabla^g_{X_h}(\phi \circ g) \\
		&= \bar \nabla_{\bar X}\bar Y^* \bbullet \phi + \bar Y^* \bbullet \bm{\bar \nabla}_{\bar X} \bar \phi.
	\end{align*}

	Finally, \cref{EqbarCliffordHermitian} follows directly from the definitions and \cref{EqCliffordHermitian}.
\end{Prf}

Now, these universal structures allow us to reformulate \cref{LemSigmagMFromSigmaMVector} as follows.

\begin{Thm}
	\label{ThmEmbDirBun}
	Let $(M, \Theta)$ be a spin manifold with fixed topological spin structure and $\bar \pi^{\Sigma}_{SM}: \bar \Sigma M \to J^1 \pi^{r,s}$ be the jetted universal spinor bundle from \cref{EqDiagUSBC}. For every metric $g$ on $M$, there exists a morphism $\bar I_g$ of vector bundles such that
	\begin{align}
		\label{EqUnivSpinorJet}
		\begin{split}
			\xymatrix{
				\Sigma^g M
					\ar@/^2pc/[drr]^-{\bar \iota_g}
					\ar@/_2pc/[ddr]_-{\pi^g_M}
					\ar@{-->}[dr]^-{\bar I_g}
				\\
				& j^1(g)^* \bar \Sigma M
					\ar[r]^-{j^1(g)^*}
					\ar[d]
				& \bar \Sigma M
					\ar[d]^{\bar \pi^{\Sigma}_{SM}}
				\\
				&M
					\ar[r]^-{j^1(g)}
				&J^1 \pi^{r,s}.
			}
		\end{split}
	\end{align}
	commutes. In addition, $\bar I_g$ is isometric with respect to the spinorial metric on $\pi^g_M$ and $j^1(g)^* \bm{\bar \eta}$, it is compatible with the Clifford multiplication $\mathfrak{m}^g$ and $j^1(g)^* \overline{\mathfrak{m}}$ and it is compatible with the spinorial Levi-Civita connection on $\pi^g_M$ and $j^1(g)^* \bm{\bar \nabla}$.
\end{Thm}

\begin{Prf}
	As in the proof of \cref{LemSigmagMFromSigmaMVector}, the commutativity of \cref{EqUnivSpinorJet} follows again from the universal property of the pullback and we obtain
	\begin{align*}
		j^1(g)^*(\bm{\bar \eta})(\bar I_g(\psi), \bar I_g(\psi'))
		& =\bm{\bar \eta}(\bar \iota_g(\psi), \bar \iota_g(\psi'))
		=\bm{\eta}(\psi, \psi')
		=g(\psi, \psi'),\\
		j^1(g)^*(\overline{\mathfrak{m}})(V \otimes \bar I_g(\psi)) &= \overline{\mathfrak{m}}((j^1(g),V) \otimes \bar \iota_g(\psi)) = V \cdot_g \psi = \mathfrak{m}^g(V \otimes \psi).
	\end{align*}
	Now, since $d\pi^{r,s}_{1,0}dj^1(g) V = d(\pi^{r,s}_{1,0} \circ j^1(g))V = dg V$, we obtain from \cref{EqDefUnivConn}
	\begin{align*}
		(j^1(g)^*\bm{\bar \nabla})_V(\bar I_g(\psi)) = \bm{\bar \nabla}_{d j^1(g) V}(\bar I_g (\psi)) = \nabla^g_V \psi,
	\end{align*}
	which implies the claim.
\end{Prf}

\subsection{Naturality questions}
\label{SubSectNaturality}

To discuss naturality aspects of the spinor bundle, it is helpful to introduce some lightweight category theoretic language and slightly rephrase what we know so far.

\begin{Def}[pseudo-Riemannian spin manifolds]
	The \emf{category of pseudo-Riemannian spin manifolds}, $\cat{pRiemSpinMf_{r,s}}$, consists of tuples $(M,g,\Theta^g)$ where $M$ is a smooth manifold, $g$ is a metric of signature $(r,s)$ and $\Theta^g$ is a metric spin structure. A \emf{morphism} between two such objects is a tuple $(f,\tilde F)$, where $f$ is an orientation-preserving isometric diffeomorphism and $\tilde F$ is a lift to the metric spin structures such that
	\begin{align}
		\label{EqDefSpinMetricMorphism}
		\begin{split}
			\xymatrix{
				\Spin^{g_1} M_1
					\ar[r]^-{\tilde F}
					\ar[d]^-{\Theta_1}
				& \Spin^{g_2} M_2
					\ar[d]^-{\Theta_2}
				\\
				\SO^{g_1} M_1
					\ar[r]^-{df}
					\ar[d]
				& \SO^{g_2} M_2
					\ar[d]
				\\
				M_1
					\ar[r]^-{f}
				& M_2
			}
		\end{split}
	\end{align}
	commutes.
\end{Def}

\begin{Def}[generalized Dirac bundles] 
	\label{DefGenDiracBdl}
	The \emf{category of (generalized) Dirac Bundles}, $\cat{DiracBdl}$, consists of tuples $(\pi_S:S \to N, h, \nabla^S, \pi_E: E \to N, g, \nabla^E, \mathfrak{m}^g)$, where $N$ is a smooth manifold, $\pi^S$ is a complex vector bundle with a metric $h$ and a metric connection $\nabla^S$, $\pi_E$ is a real vector bundle with a metric $g$ of signature $(r,s)$ and a metric connection $\nabla^E$, and $\mathfrak{m}^g:E \otimes S \to S$, $X \otimes s \mapsto X \cdot s$, is a Clifford multiplication. In addition these structures satisfy the following compatibility conditions:
	\begin{align*}
		- 2 g(X,Y) \psi &= X \cdot Y \cdot \psi + Y \cdot X \cdot \psi, \\
		\nabla_V^S(X \cdot \psi) &= \nabla_V^E X \cdot \psi + X \cdot \nabla^S_V \psi,  \\
		h(X \cdot \psi, \psi') &= (-1)^{s+1} h(\psi, X \cdot \psi'),
	\end{align*}	
	where $V \in \Gamma(\tau_N)$, $X \in \Gamma(\pi_E)$, $Y \in \Gamma(\pi_S)$. A \emf{morphism} between two generalized Dirac bundles is a tuple $(f,F, \bar F)$ such that 
	\begin{align}
		\label{EqDefMorphismDiracBundle}
		\begin{split}
			\xymatrixcolsep{3.5em}
			\xymatrix{
				S_1
					\ar[r]^-{\bar F}
					\ar[d]^{\pi_{S_1}}
				& S_2
					\ar[d]^{\pi_{S_2}}
				& \hspace{-2em}
				& E_1
					\ar[r]^-{F}
					\ar[d]^-{\pi_{E_1}}
				& E_2
					\ar[d]^-{\pi_{E_2}}
				\\
				N_1
					\ar[r]^-{f}
				& N_2,
				& \hspace{-2em}
				&N_1
					\ar[r]^-{f}
				&N_2,
			}
		\end{split}
	\end{align}
	commute. In addition, we require $f$ to be smooth, $F$ and $\bar F$ to be isometric, and $F^* \nabla^{E_2} = \nabla^{E_1}$, $\bar F^* \nabla^{S_2} = \nabla^{S_1}$, $\bar F \circ \mathfrak{m}_1 = \mathfrak{m}_2 \circ (F \otimes \bar F)$.
\end{Def}

With this notation, any Dirac bundle $S \to M$ in the classical sense, see for instance \cite{Roe}, is a generalized Dirac bundle over $E = TM$. It is well known, see for instance \cite{LM}, that associating a spinor bundle to a pseudo-Riemannian spin manifold is natural in the following sense.

\begin{Lem}
	\label{LemNaturalityIsometry}
	Let $\fun{\Sigma}:\cat{pRiemSpinMf_{r,s}} \to \cat{DiracBdl}$ be the map that assigns to each pseudo-Riemannian spin manifold $(M,g,\Theta^g)$ its classical spinor bundle $\Sigma^g M := \Spin^g \times_{\rho} \Sigma_{r,s}$ together with its Dirac structure from \cref{RemDiracStructSpingM} and its tangent bundle $\tau_M:TM \to M$. Further, we assign to any morphism $(f, \tilde F)$ of pseudo-Riemannian spin manifolds, the morphism $(f, df, \bar F)$ of Dirac bundles, where $\bar F$ is defined by setting $\bar F:\Sigma^{g_1} M_1 \to \Sigma^{g_2} M_2$, $[\tilde b, v] \mapsto [\tilde f(\tilde b), v]$. Then $\fun{\Sigma}$ is a functor.
\end{Lem}

\begin{Rem} 
	\label{codim1}
	We can extend $\fun{\Sigma}$ with little modifications to oriented isometric codimension-one immersions instead of isometric diffeomorphisms. Here, $df : \SO^{g_1} M_1 \rightarrow \SO^{g_2} M_2$ is replaced by $(df, \nu)$ which completes the pushed-forward basis by the right choice of normal vector $\nu$ to an oriented orthonormal basis, for details and applications see \cite{BaerGaudMor} or \cite{BaumMuell}. If the dimension of the hypersurface is even and we want to interpret $\fun{\Sigma}$ in a \emph{contravariant manner}, we need two spinor bundles on the hypersurface, which can be identified with the spinor bundle in the ambient space. This is due to the jump in dimension of the spinor module in this case. In the following, by a slight abuse of notation, we will use the symbol $\tilde{F}$ also for such a spin lift of an oriented isometric immersion $f$.
\end{Rem}

\begin{Rem}
	For any morphism $(f,\tilde F)$ between $(M_1, g_1, \Theta_1)$ and $(M_2, g_2, \Theta_2)$, the resulting morphism $(f, df, \bar F) = \fun{\Sigma}(f, \tilde F)$ satisfies $\Dirac^{g_2} \circ \bar F = \bar F \circ \Dirac^{g_1}$ and the two spin manifolds are Dirac isospectral. So the two spin manifolds are \emf{spin isometric} in this case. 
\end{Rem}

Now, being spin isometric is a very restrictive condition. We would like to investigate in what sense one can (not) relax this condition without losing the naturality. First, we discuss the universal spinor jet bundle from \cref{EqDiagUSBC}.

\begin{Def}[spin manifolds]
	The \emf{category of spin manifolds}, $\cat{SpinMf}$, consists of tuples $(M,\Theta)$ where $M$ is a smooth manifold and $\Theta$ is a topological spin structure. A \emf{morphism} between two such objects is a tuple $(f,\tilde F)$, where $f$ is an orientation-preserving diffeomorphism and $\tilde F$ is a lift to the topological spin structure $\GLtp M$ such that $\Theta_2 \circ \tilde F = df \circ \Theta_1$ as in \cref{EqDefSpinMetricMorphism} (with $\SO^g M$ replaced by $\GLp M$ and $\Spin^g M$ replaced by $\GLtp M$).  
\end{Def}

Recall that $\tau^{r,s}_M:(\pi^{r,s})^*(TM) \to \S_{r,s}M$ denotes the tangent bundle of $M$ pulled back to $\S_{r,s} M$. We also define $\bar \tau^{r,s}_M: (\pi^{r,s}_0)^*(TM) \to J^1\pi^{r,s}$.

\nomenclature[taubarrsM]{$\bar \tau^{r,s}_M$}{$(\pi^{r,s}_0)^*(TM) \to J^1\pi^{r,s}$}

\nomenclature[usjb]{$\fun{usjb}$}{universal spinor jet bundle functor}

\begin{Thm}
	\label{ThmUSjBisFunctor}
	The map $\fun{usjb}:\cat{SpinMf} \to \cat{DiracBdl}$, that maps a spin manifold to its universal spinor jet bundle via 
	\begin{align*}
		(M,\Theta) \mapsto (\bar \pi^{\Sigma}_{SM}:\bar \Sigma M \to J^1\pi^{r,s}, \bm{\bar \eta}, \bm{\bar \nabla}, \tau_{J^1 \pi^{r,s}}:TJ^1\pi^{r,s} \to J^1 \pi^{r,s}, \bm{\bar g}, \bar \nabla, \overline{\mathfrak{m}}),
	\end{align*}
	is a functor. In particular, for any morphism $(f, \tilde F)$ between spin manifolds $(M_j, \Theta_j)$, $j=1,2$, there exist commutative diagrams
	\begin{align}
		\label{EqUSJBnaturality}
		\begin{split}
			\xymatrixcolsep{3.5em}
			\xymatrix{
				\bar \Sigma M_1
					\ar@{-->}[r]^-{\bar F}
					\ar[d]^{\bar \pi^{\Sigma}_{SM_1}}
				& \bar \Sigma M_2
					\ar[d]^{\bar \pi^{\Sigma}_{SM_2}}
				& \hspace{-2em}
				& (\pi_0^{r,s,1})^*(TM)
					\ar@{-->}[r]^-{\hat F}
					\ar[d]^-{\bar \tau^{r,s}_{M_1}}
				& (\pi_0^{r,s,2})^*(TM)
					\ar[d]^-{\bar \tau^{r,s}_{M_1}}
				\\
				J^1 \pi^{r,s}_1
					\ar@{-->}[r]^-{\bar F_{r,s}}
				& J^1 \pi^{r,s}_2,
				& \hspace{-2em}
				&J^1 \pi^{r,s}_1
					\ar@{-->}[r]^-{\bar F_{r,s}}
				&J^1 \pi^{r,s}_2,
			}
		\end{split}
	\end{align}
	where $(\bar F_{r,s}, \bar F)$ and $(\bar F_{r,s}, \hat F)$ are isomorphisms of vector bundles compatible with the universal structures, i.e. all the relations 
	\begin{align}
		\label{EqUnivStuctsNatural}
		\bar F^*\bm{\bar \eta_2} = \bm{\bar \eta_1}, &&
		\hat F^* \bm{\bar \nabla^1} = \bm{ \bar \nabla^2}, &&
		\hat F^* \bar \nabla^1 = \bar \nabla^2, &&
		\bar F \circ \mathfrak{\overline{m}_1} = \mathfrak{\overline{m}_2} \circ (\hat F \otimes \bar F), &&
		\hat F^* \bm{\bar g_2} = \bm{\bar g_1},
	\end{align}	
	are satisfied.
\end{Thm}

\begin{Prf} 
	The pullback of tensor fields via $f^{-1}$ gives a map $F_{r,s}: \S_{r,s}M_1 \to \S_{r,s}M_2$, which induces the map $\bar F_{r,s}$ and also the map $\hat F$. The map $\Sigma M_1 \to \Sigma M_2$, $[\tilde b, v] \mapsto [\tilde F(\tilde b), v]$ is well-defined and descends to $\bar F: \bar \Sigma M_1 \to \bar \Sigma M_2$. This gives existence and commutativity of \cref{EqUSJBnaturality}. \\
	One now has to check all the compatibility conditions \cref{EqUnivStuctsNatural} in detail, which is a bit tedious. We only discuss the connection $\bm{\bar \nabla^1}$: Recall that by \cref{LemNatMetricIsom}, the map $F_{r,s}$ is an isometry. Since the vertical part of the connection is given via the metric from \cref{LemNatMetrConn}, it is preserved on the level of principal bundles. Since $\Sigma M$ is an associated bundle with induced connection, the claim follows also for $\bar \Sigma M$. The horizontal part follows from the naturality of the spinorial Levi-Civita connection with respect to spin isometries. 
\end{Prf}
	
\begin{Rem} 
	\label{DefSpinMf1}
	The modification towards oriented codimension-one immersions $f$ instead of diffeomorphisms as in \cref{codim1} works here as well, and such a lift will also be denoted by $\tilde{F}$. We denote by $\cat{SpinMf}^1$ the category that is identical to $\cat{SpinMf}$, but includes as morphisms oriented codimension-one immersions. The technicality with the jump in the dimension of the spinor module as in \cref{codim1} applies here as well.
\end{Rem}

\nomenclature[usb]{$\fun{usb}$}{universal spinor bundle functor}

\begin{Rem}[non-jetted version]
	\label{RemUSBisFunctor}
	Of course one can get a version of \cref{ThmUSjBisFunctor} for the non-jetted universal spinor bundle $\pi^{\Sigma}_{SM}: \Sigma M \to \S_{r,s} M$ and one obtains an analogous functor $\fun{usb}$. In the range category $\cat{DiracBdl}$ one would have to replace the connection by a vertical connection and one has to remove all bars in \cref{EqUnivStuctsNatural}. As a result, even the entire diagram \cref{EqDiagUSBC} is natural with respect to spin diffeomorphisms.
\end{Rem}

For the purpose of variational theory, it would be desirable to have a vector bundle \emph{over $M$} for the spinors that is natural in the category of manifolds and diffeomorphisms and on which only the secondary structure (metric, connection, Clifford multiplication) is induced by the metric of the underlying manifold. This ansatz, however, is doomed to fail as shown in the following theorem:

\begin{Thm}
	\label{ThmNoClassicalUSB}
	There is no functor $\fun{\Sigma^{top}}$ such that 
	\begin{align}
		\label{EqNatClassicSpinorBundle}
		\begin{split}
			\xymatrix{
				\cat{RiemSpinMf}
					\ar[r]^-{\fun{\Sigma}}
					\ar[d]^-{\fun{MetrTop}}
				& \cat{DiracBdl}
					\ar[d]^-{\fun{forget}}
				\\
				\cat{SpinMfds}
					\ar@{-->}[r]^-{\nexists \fun{\Sigma^{top}}}
					\ar[d]^-{\fun{forget}}
				& \cat{VB}
					\ar[dl]^-{\fun{pr}}
				\\
				\cat{Mfds}
			}
		\end{split}
	\end{align}
	is a commutative diagram of categories and functors. Here, $\fun{forget}$ denotes the functor that forgets all additional structure except the manifold and vector bundle structure, $\fun{pr}$ maps a vector bundle $E \to M$ in the category $\fun{VB}$ of vector bundles to its base $M$, and $\fun{MetrTop}$ replaces a metric spin structure by its topological one (via the fiber bundle extension corresponding to the inclusion homomorphism $\SO_m \hookrightarrow \GLp_m$).
\end{Thm}

\begin{Prf}
	Assume that \cref{EqNatClassicSpinorBundle} exists and denote by $\GLtp_m \to \GLp_m$ the non-trivial double cover. We let $\tilde b \in \GLtp_m$ be arbitrary and think of $\vartheta(\tilde b) =: b \in \GLp_m$ as an isomorphism $\R^m \to \R^m$, $x = x^i e_i \mapsto x^i b_i $. 
	Let $\sigma: \Sigma \R^m \to \R^m$ be the spinor bundle with respect to the Euclidean metric, i.e. $\sigma := \fun{\Sigma^{top}}(\fun{MetrTop}(\R^m,\Spin_m,\bar g)) = \fun{\Sigma^{top}}(\R^m, \GLtp_m)$. We set $V := \sigma^{-1} (0)$ and obtain an isomorphism $\tilde r_b :\fun{\Sigma^{top}}(b, \tilde b)|_V:V \to V$. All in all, we obtain a finite dimensional representation $\tilde r:\GLtp_m \to \GL(V)$, $\tilde b \mapsto \tilde r_b$. By \cite[Lem. II.5.23]{LM}, this representation descends to a representation of $\GLp_m$, i.e. there exists a representation $r:\GLp_m \to \GL(V)$ such that
	\begin{align}
		\label{EqReprFactors}
		\begin{split}
			\xymatrix{
				\GLtp_m
					\ar[r]^-{\tilde r}
					\ar[d]^-{\vartheta_m}_-{2:1}
				& \GL(V)
				\\
				\GLp_m
					\ar@{-->}[ur]_-{r}
			}
		\end{split}
	\end{align}
	commutes.
	Now, consider the unit element $e \in \GLp_m$ and obtain $\vartheta_m^{-1}(1) =: \{ \tilde e_+, \tilde e_- \}$, where $\tilde e_+$ is the unit element in $\GLtp_m$. Since \cref{EqReprFactors} commutes, we obtain that 
	\begin{align*}
		f:= \tilde r(\tilde e_{-}) = (r \circ \vartheta_m)(\tilde e_{-}) = r(\tilde e_+) = \id_V.
	\end{align*}
	Notice that $\tilde e_{-} \in \Spin_m$. Therefore, we can consider $\fun{forget}(\fun{\Sigma}(\id_{\R^m}, \tilde e_{-}))|_V = f \in \GL(V)$. 
	By construction $\Sigma^{\bar g} \R^m = \Spin_m \times_{\rho} \Sigma_m$, where $\rho$ is the standard spin representation, and the morphism $f$ is explicitly given by $\Sigma \R^m \to \Sigma \R^m$, $[\tilde b, v] \mapsto [\tilde b.\tilde e_{-}, v] = [\tilde b, \rho(\tilde e_-)^{-1} v]$. Thus $f \neq \id_V$, since $\rho$ is faithful, which is a contradiction.
\end{Prf}

\subsection{Universal Dirac operator}

The naturality of the universal spinor bundle induces a naturality assertion for the Dirac operator. Recall that if $\pi_P: P \to M$ and $\pi_Q:Q \to M$ are fibre bundles, then a map $D:\Gamma_{\loc}(\pi_P) \to \Gamma_{\loc}(\pi_Q)$ is a \emf{differential operator of order $k$}, if there exists a bundle morphism $\sigma(D):J^k \pi_P \to Q$ over $\id_M$, called \emf{symbol of $D$}, such that $D(s) = \sigma(D)(j^k(s))$ for any $s \in \Gamma_{\loc}(\pi_P)$, see for instance \cite[Def. 6.2.22]{saunders}. 

\begin{Thm}[universal Dirac operator]
	\label{ThmDiracNaturalDiffOp} $ $ 
	\begin{enumerate}
		\item
			There exists a natural linear differential operator $\bar \Dirac$ of order $1$ on $\bar\pi^\Sigma_{SM}$, which for every metric $g$ pulls back to the usual Dirac operator, i.e. there exists a commutative diagram of vector bundle homomorphisms \footnote{As in our case of double bundles, the projections of the respective bundles are not uniquely given by the total spaces, we use the not so common but very practical way of writing a diagram  of vector bundle homomorphims as a diagram {\em between the projections}, meaning that there is a corresponding commuting diagram between the respective total spaces preserving the and linear on the respective fibers. }
			\begin{align}
				\label{EqUniDiracOp}
				\begin{split}
					\xymatrixcolsep{6em}
						\xymatrix{
							j_{1,0}^1 \pi^g_M
								\ar[r]^{j^1 \bar{\iota}^g}
								\ar[d]^{\sigma(\nabla^g)}
								\ar@/_6em/[dd]_-{\sigma(\Dirac^g)}
							& j_{1,0}^1 \bar \pi^\Sigma_{SM} 
								\ar[d]^-{\sigma(\bar{\bm \nabla})}
								\ar@/^6em/@{-->}[dd]^-{\sigma(\bar \Dirac)}
							\\
							\tau^*_{M} \otimes \pi^g_M
								\ar[r]^-{(d \pi_0)^* \otimes \bar{\iota}^g}
								\ar[d]^{\mathfrak{m}^g \circ \sharp_g}
							& \tau^*_{J^1 \pi^{r,s}} \otimes \bar \pi^\Sigma_{SM} 
								\ar[d]^{a}
							\\
							\pi^g_M
								\ar[r]^-{ \bar \iota_g}
							& \bar \pi^\Sigma_{SM} 
						}
				\end{split}
			\end{align}
		
		\item
			We can also interpret $\bar{\Dirac}$ as differential operator $\Dirac$ on $\pi^\Sigma_M$ in a natural way, given explicitly as $\Dirac: \Gamma(\pi^{\Sigma}_{M}) \to \Gamma(\pi^{\Sigma}_{M})$, $\Phi=(g,\varphi) \mapsto (g, \Dirac^g \varphi)$.
	\end{enumerate}
\end{Thm}

\begin{Prf} $ $
	\begin{enumerate}
		\item
			The left hand side is the usual Dirac operator constructed via the usual Dirac structure, c.f. \cref{RemDiracStructSpingM}. The objects on the right hand side are defined as in \cref{EqDiagUSBC} and the map $\bar \iota_g$ is given by \cref{EqUnivSpinorJet}. It remains only to construct the map $a$: Let $h:(\pi^{r,s})^*(TM) \to TJ^1\pi^{r,s}$, $(j^1_x(g), V) \mapsto dg|_x(V)$, be the holonomic lift, c.f. \cref{EqJetBundleDecomposition}. This induces a map $h^*:T^*J^1\pi^{r,s} \to ((\pi^{r,s})^*(TM))^*$. The map $(\pi_{1,0}^{r,s})^*: (\pi^{r,s}_{0})^*(TM) \to (\pi^{r,s})^*(TM)$ yields a dual map $((\pi_{1,0}^{r,s})^*)^*: ((\pi^{r,s})^*(TM))^* \to ((\pi^{r,s}_{0})^*(TM))^*$. Using the universal metric from \cref{EqDefUniMetricgbar}, we obtain an isomorphism $\sharp_{\mathbf{\bar g}}:((\pi^{r,s}_0)^*(TM))^* \to (\pi^{r,s}_0)^*(TM)$. By definition, $\overline{\mathfrak{m}}:(\pi_0^{r,s})^*(TM) \otimes \bar \Sigma M \to \bar \Sigma M$. Therefore, we can set $a:= \overline{\mathfrak{m}} \circ \sharp_{\mathbf{\bar g}} \circ ((\pi^{r,s}_{1,0})^*)^* \circ h^*$. The diagram commutes by construction, in particular \cref{EqDefUnivConn}. This shows that $\sigma(\bar \Dirac)$ is natural and thus $\bar \Dirac := \sigma(\bar \Dirac) \circ j^1$ is a natural linear differential operator of order $1$ on $\bar\pi^\Sigma_SM$.  
		\item
			We can also interpret it as a differential operator on the fiber bundle $\pi^\Sigma_M $ by taking into account that there is a natural fiber bundle isomorphism identifying $ \pi_0^{r,s} \circ j^1_{0} \bar \pi^\Sigma_{SM} $ with $j^1_0 \pi^\Sigma_M$. If we recall from \cref{LemUniversalSpinorFields} that a universal spinor field $\Phi$ can be thought of as a tuple $(g, \varphi)$, this operator can also be given explicitly by $\Dirac: \Gamma(\pi^{\Sigma}_{M}) \to \Gamma(\pi^{\Sigma}_{M})$, $\Phi=(g,\varphi) \mapsto (g, \Dirac^g \varphi)$, is a differential operator of order $1$ that is natural with respect to spin diffeomorphisms as well. This can also be seen using its definition, the coordinate formula \cref{EqDiracghLocal} below and \cref{ThmUSjBisFunctor} (or rather its non-jetted version, see \cref{RemUSBisFunctor}).
	\end{enumerate}
\end{Prf}

\section{The Einstein-Dirac-Maxwell equation}
\label{SectEinsteinDirac}

Using the universal spinor bundle $\pi^{\Sigma}_M: \Sigma M \to M $, one can now formulate Einstein-Dirac-Maxwell theory as a variational problem on sections of a finite dimensional fiber bundle in any signature $(r,s)$. The advantage of this approach is that the spaces of sections $W^{k,p}(\pi^{\Sigma}_M)$ carry canonical topologies as Fréchet manifolds (Banach if $k < \infty$). At the first sight, including Maxwell fields might seem to be an unnecessary complication. However, a good reason for doing so is that, whereas solutions of Einstein-Maxwell theory can be made solutions of Einstein-Dirac-Maxwell theory by including zero spinor fields, the same is not true for Einstein-Dirac solutions: The presence of nonzero spinors entails in general the presence of nonzero Maxwell fields.

\subsection{Notation and basic definitions}

To formulate the Einstein-Dirac-Maxwell equation, we introduce the following notions. Let $g \in \mathcal{S}_{r,s}(M)$ be a metric, $\psi \in \Gamma(\pi^g_M)$ be a spinor field and $A \in \Omega^1(M)$ (thought of as a connection form on a trivial $U(1)$-bundle over $M$). Each $1$-form $A$ induces a connection $\nabla^{g,A}$ on the spinor bundle via $\nabla_X^{g,A} \psi := \nabla_X^g \psi + i A(X) \cdot \psi$. The resulting Dirac operator reads $\Dirac^{g,A} \psi = \Dirac^g \psi - A \cdot \psi$ (using Clifford multiplication on the cotangent bundle). For $q \in \R$, we define the tensor fields $T^1_{(q,g,\psi,A)}, T^2_{(g,A)} \in \Gamma(\tau^2_M)$ by
\begin{align*}
	T^1_{(q, g,\psi, A)}(X,Y) & := \tfrac{1}{2} \Re \langle X \cdot \nabla_Y^{g, q A} \psi + Y \cdot \nabla_X^{g, q A} \psi, \psi \rangle ,  \\
	T_{(g, A)}^2 &:= \tr_{(1,3)}^g(F^A \otimes F^A) - \tfrac{1}{4} g(F^A, F^A)  g, 
\end{align*}
for all $X, Y \in \Gamma(\tau_M)$. Here $F^A := dA$, which is the curvature of $A$. The metric $g$ is extended canonically to all tensor powers of the tangent bundle $\tau_M: TM \to M$. We also set $T_{(q, g,\psi,A)} := T^1_{(q, g,\psi, A)} + T^2_{(g, A)}$, the \emf{total energy-momentum tensor of $(q, g,\psi,A)$}. The field $j_\psi \in \Omega^1(M)$, defined by	$j_\psi (X) := \langle X \cdot \psi, \psi \rangle$, is called the \emf{Dirac current of $\psi$}. 

\nomenclature[Tgpsi]{$T_{(g,\psi)}$}{energy momentum tensor of $(g,\psi)$}
\nomenclature[TgA]{$T_{(g,A)}$}{energy momentum tensor of $(g,A)$}
\nomenclature[TgpsiA]{$T_{(g,\psi,A)}$}{total energy momentum tensor of $(g,\psi,A)$}
\nomenclature[jpsi]{$j_\psi$}{Dirac current}

\begin{Def}[Einstein-Dirac-Maxwell equation] 
	Let $(M,\Theta)$ be a spin manifold. For any $\lambda, q \in \R$, the system of equations 
	\begin{align}
		\label{EqEinsteinDirac}
		\begin{split}
			\Dirac^{g, q A} \psi & = \lambda \psi, \\ 
			\Ric^{g} - \tfrac{1}{2} \scal^{g} g &=  T_{(g,\psi,A)}, \\
			\delta^g F^A &= q j_\psi,
		\end{split}
	\end{align}
	is called \emf{$(\lambda, q)$-Einstein-Dirac-Maxwell equation}. $q$ is called the \emf{charge} and $\lambda$ is called the \emf{mass}. We want to solve this equation for $g$, $\psi$, $A$, and denote such a solution by $(M, g, \psi, A)$.
\end{Def}

\begin{Rem}[neutral systems]
	\label{RemNeutralSystems}
	Physically reasonable systems often consist of various particles of different charges and masses. By summation of the Dirac currents in the last equation and of energy-momentum tensors $T^1_{(q_i, g, \psi_i, A)}$ in the second equation, the equations can readily be generalized to a Dirac system comprising various spinor fields, for details cf. \cite{GM}. This is important as it allows to consider \emf{neutral systems}, which are systems with $\sum_i q_i =0$. If, moreover, $\lambda_i = 0$ for all $i$, these systems have a well-posed initial value formulation for small initial values in asymptotically flat spacetimes, cf. \cite{GM}.
\end{Rem}

\begin{Def}[universal Einstein-Dirac-Maxwell operator]
	The map $ \EDM_{(\lambda, q)} : \Gamma(\pi^{\Sigma}_M \oplus \tau^*_M) \to \Gamma(\tau^*_M \otimes \tau^*_M \oplus \pi^{\Sigma}_M \oplus \tau_M^*)$,
	\begin{align}
		\label{EqDefEinsteinDiracMaxwellOp}
		\Phi=(g,\psi,A) \mapsto (\Ric^g - \tfrac{1}{2} \scal^g g - T_{(q,g,\psi,A)} , 
		(g, \Dirac^{g, q A} \psi - \lambda \psi), 
		\delta^g F^A - q j_\psi),
	\end{align}
	is called \emf{universal Einstein-Dirac-Maxwell operator}. Obvious extensions of the above operators to $N$ spinor fields are also denoted by $\EDM_{(\lambda, q)}$; here $\lambda = (\lambda_1, \ldots , \lambda_N )$, $q= (q_1, \ldots , q_N)$. 
\end{Def}

\begin{Thm}
	\label{ThmEDMNaturalDiffOp}
	The operator $\EDM_{(\lambda, q)} $ is a natural differential operator whose zeroes correspond exactly to the solutions of the Einstein-Dirac-Maxwell equations.
\end{Thm}

\begin{Prf}
	In the case of $N$ spinor fields, one has to replace $\pi^{\Sigma}_{SM}$ from the universal spinor bundle in \cref{EqUniSpinDouble} by its $k$-fold fiber product. For simplicity, we consider only the case $N=1$. Then \cref{ThmDiracNaturalDiffOp} gives the result for the Dirac-part of $\EDM$ and the statements for the other parts are standard.
\end{Prf}

\subsection{Solutions as critical points}

Now, want to express solutions of \cref{EqEinsteinDirac} as critical points of a functional to be defined in \cref{DefEinsteinDiracMaxwellFunctional}. To that end we require certain trivializations, which can be obtained from identification isomorphisms, which are an important application of the partial connection on the universal spinor bundle anyway. 

\begin{Def}[identification isomorphisms]
	\label{Defbetagh}
	Let $\mathbf{g}: I \to \mathcal{S}_{r,s}(M)$, $t \mapsto g_t$, be a smooth path of metrics, $g:=g_0$, $h:=g_1$. For any $x \in M$ and any $\Phi \in \Sigma M|_{g(x)}$, let $\beta_{g,h}(\Phi) \in \Sigma M|_{h(x)}$ be obtained by $\bm{\nabla}$-parallel transport of $\Phi$ along $\mathbf{g}$. The resulting map $\beta_{g,h}:\Sigma^{g} M \to \Sigma^{h} M$ is called an \emf{identification isomorphism}.
\end{Def}

\nomenclature[betagh]{$\beta_{g,h}$}{identification isomorphism}

These identification maps satisfy $\beta_{g,h} \circ \beta_{h,g} = \id$, and one can use them to pull back any Dirac operator $\Dirac^h$ to $\Gamma(\pi^g_M)$. This operator can be expressed more explicitly as follows: For any metric $g \in \mathcal{S}_{r,s}(M)$, we denote by $\flat_g:TM \to T^*M$ the musical isomorphism and by $\sharp_g:T^*M \to TM$ its inverse. If $h \in \mathcal{S}_{r,s}(M)$ is any other metric, the isomorphism $a_{g,h} := \sharp_g \circ \flat_h \in \Iso M$ satisfies
\begin{align}
	\label{EqDefagh}
	\forall X, Y \in TM: g(a_{g,h}(X), Y) = h(X,Y).
\end{align}
Clearly, if $g=h$ we obtain $a_{g,h} = \id$, which is positive definite. Therefore, if $h$ is in a small neighborhood of $g$, the map $a_{g,h}$ is still positive definite. (In case $g,h$ are Riemannian, $a_{g,h}$ is always positive definite.) The map $b_{g,h} := \sqrt{a_{g,h}^{-1}}$ (in the sense of a positive definite square root) satisfies
\begin{align}
	\label{EqDefbgh}
	\forall X, Y \in TM: h(b_{g,h}(X),b_{g,h}(Y)) = g (X,Y),
\end{align}
thus it maps a $g$-pseudo-orthonormal basis to an $h$-pseudo-orthonormal basis. 

Any identification isomorphism $\beta_{g,h}$ from \cref{Defbetagh} induces a map on sections $\Gamma(\pi^g_M) \to \Gamma(\pi^h_M)$ also denoted by $\beta_{g,h}$. The problem is that for any two metrics $g,h \in \mathcal{S}_{r,s}(M)$ there might be no path in $\mathcal{S}_{r,s}(M)$ joining them. Even if they can be joined, the path is not unique. To handle this problem, we introduce the following notion.

\begin{Def}
	Let $M$ be a (possibly non-compact) manifold and $g,h \in \mathcal{S}_{r,s}(M)$. We say $g$ and $h$ are \emf{joinable}, if the path $g_t := g + t(h-g)$, $t \in [0,1]$, is contained in $\mathcal{S}_{r,s}(M)$ and for any $t \in [0,1]$, the map $a_{g,g_t}$ is positive definite.
\end{Def}

Obviously, if we are given a compact subset $C \subset M$, a pseudo-Riemannian metric $h$ and an auxiliary Riemannian metric $k$, there is a positive number $a_C$ such that all pseudo-Riemannian metrics $g$ coinciding with $h$ outside of $C$ and with $\|g-h \|_k < a_C$ are joinable to $h$. Even more, there is a positive number $b_C$ such that all pseudo-Riemannian metrics $g$ coinciding with $h$ outside of $C$ and with $\|g-h \|_k < a_C$ are joinable to each other. We will call such a neighborhood {\bf $C$-convex}.
For any two joinable $g,h$, we obtain a unique identification isomorphism $\beta_{g,h}$ from \cref{Defbetagh}. 

\begin{Thm}
	\label{ThmDiracghLocal}
	Let $g,h \in \mathcal{S}_{r,s}(M)$ be joinable. Then the operator $\Dirac^{h}_{g} := \beta_{h,g} \circ \Dirac^h \circ \beta_{g,h}: \Gamma(\pi^g_M) \to \Gamma(\pi^g_M)$ has the local coordinate representation
	\begin{align} 
		\label{EqDiracghLocal}
		\begin{split}
			\Dirac^{h}_{g} \psi
			& = \sum_{i=1}^{m}{e_i \cdot \nabla^g_{b_{g,h}e_i}\psi}
			+\frac{1}{4} \sum_{i,j=1}^{m}{e_i \cdot e_j \cdot \left( b_{h,g} (\nabla^h_{b_{g,h}(e_i)}(b_{g,h}e_j)) - \nabla^g_{b_{g,h}(e_i)}e_j \right) \cdot \psi},
		\end{split}
	\end{align}
	where $e_1, \ldots, e_m$ is a local pseudo-orthonormal frame and $\psi \in \Gamma(\pi^g_M)$ is any spinor field.
\end{Thm}

\begin{Prf}
	In the case of a compact Riemannian spin manifold, this formula is exactly \cite[Thm. 20]{BourgGaud} and holds for any two Riemannian metrics $g$ and $h$. The proof goes through in the general case, since we required $g$ and $h$ to be joinable. 
\end{Prf}

For the rest of this subsection, we assume the manifold $M$ to be compact.

In \cite{FriedEinstDirac}, using the variation formula for the Dirac operator established by Bourguignon and Gauduchon, Kim and Friedrich show how solutions of the Einstein-Dirac equation in the Riemannian case can be identified with \emph{stationary} points of some functional w.r.t. special variations. Now we show how to characterize any solution $(M, g,\psi,A)$ of the Einstein-Dirac-Maxwell equation as a \emph{critical} point, i.e. in the sense of $d_{(M, g,\psi,A)} \EDM=0$, of a functional $\EDM$ defined on the sections of the finite dimensional fiber bundle $\cat{usjb}$ from \cref{EqUniSpinDouble}. We denote by $\Lambda M$ the bundle of exterior forms on $M$.

\begin{Def}[Einstein-Dirac-Maxwell functional]
	\label{DefEinsteinDiracMaxwellFunctional}
	The \emf{$(q, \lambda)$-Einstein-Dirac functional} is given by$L_{(q, \lambda)}:\Gamma(\pi^{\Sigma}_M \oplus \Lambda^1 M) \to \R$, $(\Phi=(g,\psi), A) \mapsto \int_{M}{\mathcal{L}_{(q, \lambda)} \circ j^2(\Phi, A)}$, where $\mathcal{L}_{(q, \lambda)} \in \mathcal{C}^{\infty}(J^2(\pi^{\Sigma}_M \oplus \Lambda^1 M), \Lambda^m M)$ is defined by 
	\begin{align*}
		\mathcal{L}_{(q, \lambda)}(j^2(\Phi, A)) 
		:= \Big{(} \scal^{g} 
		+ \lambda \langle \psi, \psi \rangle 
		- \langle \Dirac^{g, qA} \psi, \psi \rangle
		- \frac{1}{2} g(F^A, F^A) \Big{)} dv^g.
	\end{align*}
	For simplicity, we set $L := L_{(\lambda, q)}$.
\end{Def}

\begin{Rem}[characterization of critical points]
	\label{RemCharCritical}
	Let $\pi:P \to M$ be a fiber bundle, $\tau_P^v:T^v P \to M$ be the vertical tangent bundle, and $L:\Gamma(\pi) \to \R$ be a functional. A section $s \in \Gamma(\pi)$ is a critical point of $L$ if and only if $\tfrac{d}{dt} s_t |_{t=0} = 0$ for  any $1$-parameter family $s_t \in \Gamma(\pi)$ that is smooth in $t$. Using the identification between $\Gamma(s^* \tau_P^v)$ and $T_s \Gamma(\pi)$, we get that $s$ is critical if and only if $\tfrac{d}{dt} Fl_{X}^t(s) |_{t=0} = 0$ for all $X \in \Gamma(\tau_{\Gamma(\pi)})$. Here, $Fl_{X}$ denotes the flow of $X$. In case that $L = \mathcal{L} \circ j^2$ as in \cref{DefEinsteinDiracMaxwellFunctional}, this can be written equivalently as $\int_M{(j^2 X)(\mathcal{L})} = 0$.
\end{Rem}

\begin{Rem}
	\label{RemProducStructKF}
	Now, we want to examine critical points and we want to justify the product structure as in \cite{FriedEinstDirac}: As we assumed $M$ to be compact, there is an open covering of $S_{r,s} M$ by $M$-convex, that is, convex, subsets. In such a subset $U$ we trivialize the bundle $\pi^\Sigma_{SM}$ by means of the parallel transport to a fixed metric $ g \in U $: For any $e \in \Sigma M$, we define $ t_g (e) :=  ( \pi^\Sigma_{SM} (e) , \beta_{\pi^\Sigma_{SM} (e) , g} I_g^{-1} (e)) $. Consider the open subset set $W$ of $\Gamma(\pi^\Sigma_M)$ consisting of all the sections $\Psi $ of $\pi^\Sigma_M$ such that $\pi^\Sigma_S  \circ \Psi \in U $. Then $K_g := t_g \circ \_$ is a Fr\'echet diffeomorphism between $W$ and $ U \times \Gamma (\Sigma^g M)$. If $g,h \in U$, then $K_h \circ K_g^{-1} (u, \psi) = (u, \beta_{u,h} \circ \beta_{g,u} (\psi) )$ for all $u \in U, \psi \in \Gamma(\Sigma^u M)$. Consequently, if $\Psi$ is a critical point of $\EDM \circ K_g$, then $\Psi$ is a critical point of $\EDM \circ K_{\pi^\Sigma_S \Psi}$, so we can show criticality in a trivialization w.r.t. $\pi^\Sigma_S \Psi$, as it is done in \cite{FriedEinstDirac}.
\end{Rem}

\begin{Thm}
	\label{ThmCritEinstDiracFried}
	Let $M$ be compact. A tuple $(\Phi, A)$, $\Phi=(g,\psi)$, is a critical point of $L$ if and only if $(M, g, \psi, A)$ is a solution of the Einstein-Dirac-Maxwell equation.
\end{Thm}

\begin{Prf}
	Since the inclusion of an additional Maxwell field is trivial, it suffices to discuss the Einstein-Dirac equation (which corresponds to $A=0$).
	First of all, note that \cite[Thm. 2.1]{FriedEinstDirac} holds for any signature of $g$, since the proof goes through without modifications. The theorem states that a spinor field $(g, \psi, 0)$ satisfies the Einstein-Dirac equation if and only if the functional $L(g, \psi, 0)$ is stationary under all variations of type 	
	\begin{align}
		\label{EqVarKimFried}
		t \mapsto (g + tk, \psi, 0), &&
		t \mapsto (g, \psi + t \varphi,0 ),
	\end{align}
	where the first variation is to be understood using the identifications $\beta_{g,h}$. Therefore, we can identify these variations with $t \mapsto \Phi_t := (g + tk, \beta_{g,g +tk} \varphi)$ in our setting, see \cref{RemProducStructKF} The second variation in \eqref{EqVarKimFried} is already a variation of a universal spinor field in the spinor direction. \\
	Now let us show the ``$\Longrightarrow $'' direction: Assume that $(\Phi, 0) = (g, \psi, 0)$ is a critical point of $L$. In particular, it is critical under the variation \cref{EqVarKimFried}. Therefore, using \cite[Thm. 2.1]{FriedEinstDirac} it is a solution of the Einstein-Dirac equation. \\	
	Now, for the other direction assume that $\Phi = (g, \psi)$ is a solution of the Einstein-Dirac equation. Then \cite[Thm. 2.1]{FriedEinstDirac} implies that it is critical under variations of \cref{EqVarKimFried}. We have to show that those generate all the variations by curves in $\Gamma(\pi^{\Sigma}_M)$. By \cref{RemCharCritical} these correspond to $\pi^{\Sigma}_{M}$-vertical vector fields on $\Sigma M$.	
	We need that $dL(\Phi) X^v = 0$ for all $X^v \in \Phi^*(\tau_{\Sigma M}^v)$, where $\tau^v_{\Sigma M} := \tau_{\Sigma M}|_{(d \pi^{\Sigma})^{-1}(0)}$. The Bourguignon-Gauduchon connection gives us a global decomposition $\tau_{\Sigma M}^v  = \tau_{\Sigma M}^{vv} \oplus \tau_{\Sigma M}^{vh}$ and corresponding projections $\pi^{vv}$ and $\pi^{vh}$. Now $X^{vv} := \pi^{vv}(X^v)$ is the variational vector field of a variation of second type of \cref{EqVarKimFried} and $X^{vh} := \pi^{vh}(X^v)$ is the variational vector field of a variation of first type of \cref{EqVarKimFried}. By linearity of the derivative and by $d_{\Phi}L (X^{vv}) = 0$ and $d_{\Phi} L(X^{vh}) = 0$, we obtain $d_{\Phi} L (X^v) = 0$. 
\end{Prf}

\begin{Rem}
	The description of the solution as critical points of a Lagrangian functional on a finite-dimensional jet bundle now enables us to use all the available tools in this setting. For instance, if one wants to solve extension problems one can ask for Palais-Smale condition or if the functional is a Morse function in some appropriate sense and then use tools like the Mountain Pass Lemma.
\end{Rem}

\begin{Rem}[inclusion of boundary values]
	\label{RemCharCriticalInitial}
	Now, assume that $M$ is a compact manifold with boundary. Then we can repeat the proof of \cref{ThmCritEinstDiracFried} and show criticality in the space $\Gamma_Z(\pi^{\Sigma}_M\oplus\Lambda^1)$ of fixed boundary values $Z$.
\end{Rem}

\subsection{Existence of a maximal Cauchy development}

From now on, we focus on the Lorentzian case solely. The question of a maximal Cauchy development for Einstein-Maxwell theory has been positively answered by reducing the system \cref{EqEinsteinDirac} to a symmetric-hyperbolic operator, see for instance \cite{CB}. Now we do the same for Einstein-Dirac-Maxwell Theory. In \cite{FR}, it is indicated how to do so in the framework of two-spinor calculus. As far as we can see, this framework requires the choice of an (unnatural) trivialization of the spinor bundle and is moreover restricted to dimension $4$. For the definition of a maximal Cauchy development, we need the stronger naturality properties of \cref{LemNaturalityIsometry} to prove geometric uniqueness. Therefore, we choose another approach.

\begin{Def}[initial values]
	Let $Z = (S, g_0, K, \psi_0, A_0, A_1 )$ be a tuple, where $S$ is a $(m-1)$-dimensional spin manifold, $g_0$ is a Riemannian metric on $S$, $K \in \Gamma(\tau^2_S)$ is a symmetric $2$-tensor on $S$, $\psi_0 \in \Gamma(\pi^g_S)$ is a spinor field (for $m$ even) or two spinor fields (for $m$ odd), $A_0, A_1 \in \Gamma(\tau^*_{S \times \R}|^{S \times \{0\}})$ representing\footnote{Here, for a bundle $\pi: E \rightarrow M$ and a subset $A \subset M$, we employ the definition $\pi \vert^A := \pi \vert_{\pi^{-1} (A)}$.} initial values for the Maxwell field. Then $Z$ is called an \emf{initial value} for the $(\lambda, q)$-Einstein-Dirac-Maxwell equation if and only if the initial constraint equations\footnote{These equations are apparently not yet formulated in the first jet bundle of a geometric bundle on $S $ itself but can, by a well-known procedure, easily be reformulated as equations in the spinor bundle on $S$ and the first jet bundle of Maxwell and metric fields on $S$, for details see for instance \cite{CB} and \cite{GM}.}
	\begin{align*}
		\scal^{S,g} + (\tr_g K)^2 - \vert K \vert_g
		& = 16 \pi \; T_{(q, g, \psi, A)} (\nu,\nu), \\
		\tr_{12} (\nabla K ) - \tr_{23} ( \nabla K  ) &= 8 \pi \;  T_{(q, g, \psi, A)} (\nu,\cdot), 
	\end{align*}
	are satisfied, where $\nu$ is a future-directed unit normal vector field for $S$.
\end{Def}

\begin{Rem}
	It is well-known, see for instance \cite{CB}, that restricting a section that is a solution to the Einstein-Dirac-Maxwell equation to a Cauchy surface yields an initial value (where spinors are re-interpreted as spinors of the submanifold in the sense of \cref{codim1}). 
\end{Rem}

\begin{Def}[Cauchy development]
	\label{DefCauchyDev}
	Let $Z = (S, g_0, K, \psi_0, A_0, A_1 )$ be an initial value. Then a \emf{Cauchy development} of $Z$ is a tuple $(M, g, \psi, A, f)$ such that $(M, g, \psi, A)$ is a solution of the Einstein-Dirac-Maxwell equation, $(M,g)$ is a globally hyperbolic spin $m$-manifold and $f$ is a spin diffeomorphism $S \to f(S) \subset M$ such that 
	\begin{enumerate}
		\item 
			\label{DefCauchyDevItCauchy}
			$f(S)$ is a Cauchy surface of $(M,g)$,
		\item
			\label{DefCauchyItPullback}
			$f^*g = g_0$, $f^*W= K$ and $f^*\psi = \psi_0$, $f^* A= A_0$, $f^* \nabla_\nu^g A = A_1 $,
	\end{enumerate}
	Here, $W$ denotes the second fundamental form of $S$ (``Weingarten tensor'') and $\nu$ a normal vector field. A development is called \emf{maximal}, if for every tuple $(M',g', \psi', A', f')$ satisfying items \ref{DefCauchyDevItCauchy} and \ref{DefCauchyItPullback}, we have a spin diffeomorphism $\alpha: M' \rightarrow \alpha(M') \subset M$ with $\alpha^*(g) = g' $, $\tilde \alpha(\psi') = \psi$ (in the sense of \cref{EqDefSpinMetricMorphism}), $\alpha^* A = A'$ and $d(\alpha \circ f') = d f$.
\end{Def}

In order to discuss the existence and uniqueness of solutions of \cref{EqEinsteinDirac}, we require some notions from the theory of symmetric hyperbolic quasilinear systems, see \cite{GM} for a definition. In particular, we choose a \emph{gauge}. 

\begin{Def}[wave-gauged]
	For two Lorentzian manifolds $(N,h)$ and $(M,g)$, a map $f: N \to M$ is called a \emf{wave map}, if $\tr_h (\nabla df) = 0$. Here we consider $df \in \Gamma(\tau^*_N \otimes f^* \tau_M)$, and $\nabla$ is the natural connection on $\tau^*_N \otimes f^* \tau_M$. Let $h$ be a Lorentzian metric on an open subset $U'$ of $M$ with Cauchy hypersurface $S$. A metric $h$ on a neighborhood $U \subset U'$ of $S$ is called \emf{wave-gauged to $g$ (on $U$)}, if the identity on $U$ is a wave map from $(U, h|_U)$ to $(U,g)$. 
\end{Def}

If $(M, g)$ is an open set of the Minkowski space $\R^{1, m}$, then $h$ is in $g$-wave gauge if and only every coordinate is a $h$-harmonic function. Locally every metric can be brought into such a $g$-wave-gauge, c.f. \cite[VI.7.4]{CG}.

\begin{Lem}
\label{EverythingIsHarmonic}
	Let $(N,h)$ be a globally hyperbolic manifold containing a Cauchy hypersurface $S$. Let $(M,g)$ be a Lorentzian manifold and $f:S \to M$ be an isometric embedding such that $f(S)$ is a Cauchy surface of $M$. Then there is a neighborhood $U$ of $S$ and a diffeomorphism $\alpha:U \to \alpha(U) \subset M$ with $d\alpha |_{TS} = df|_{TS}$ such that $h$ is $\alpha^*g$-wave-gauged. For every globally hyperbolic $U$ the diffeomorphism is unique and has a lift to the spin bundle that is uniquely given by its restriction to $S$.
\end{Lem}

\begin{Prf}
	The equation for $h$ to be $\alpha^*g$-wave-gauged is $Q^h(g) := \tr_h (\nabla^g- \nabla^h) = 0$. This is a quasilinear second degree equation for $g$ and thus locally solvable for first order initial values by the usual arguments via the symmetric-hyperbolic first order prolongation of $Q^h$ and then using, e.g. \cite[16.1.2-16.1.7 and 16.2.1]{taylor3} for smooth coefficients, or \cite[Thm. 4.2]{GM} for coefficients of finite regularity. Here, we need the same statement as in the references above, but for a range $G \subset M \times \R^N$ instead of $M \times \R^N$, such that for all $p \in M$, the set $G \cap (\{p\} \times \R^N)$ is convex and reflection symmetric around $0$. By a well-known procedure \cite{CB}, we can extend the coefficients $A_j$ and $B$ as in \cite{taylor3, GM} outside $G$ by an appropriate constant map and use stability. $G$ is then chosen in such a way that for any $g$ with $g_p \in G_p$ for all $p$ satisfies that the coefficient $A_0$ as in \cite{GM} is uniformly positive. The diffeomorphism $\alpha$ is spin, i.e. has a lift to the spin bundle as the identity on $S$ has and because of the homotopy lifting property for the spin bundle, as $U$ is diffeomorphic to $S \times \R$. 
\end{Prf}

\begin{Def}[reduced Einstein-Dirac-Maxwell operator]
	Let $h$ be a globally hyperbolic metric on a manifold such that $S$ is a Cauchy hypersurface. Let $U$ be a neighborhood of $h$ consisting of metrics joinable to $h$. Then
	\begin{align*}
		\EDM_{(q, \lambda)}^h : & U \oplus \Gamma(\Sigma^h M) \oplus \Omega^1(M) \to \tau^2_M \oplus \Gamma(\Sigma^h M) \oplus \Omega^1(M)	 \\
		(g, \psi^h, A) & \mapsto
		(\Ric^g - \tfrac{1}{2} \scal^g g - T_{(g,\beta_{h,g}(\psi),q A)}
		, \square_g A - q j_{\psi}
		,  \beta_{g, h} (i \nu \cdot_g ( \Dirac^{g,q A} - \lambda)  \beta_{h, g} \psi^{h}) )
	\end{align*}
	is called \emf{reduced Einstein-Dirac-Maxwell operator}. 
\end{Def}
 
We want to show that $\EDM^h_{(q, \lambda)}$ is a symmetric hyperbolic quasilinear operator. 
 
By \cref{EverythingIsHarmonic} and by bijectivity of $\beta_{g,h}$ and the fact that every tuple $(A, \psi)$ can be brought into Lorenz gauge replacing it by $(A + df, e^{if} \psi)$ for some function $f$, we obtain a one-to-one correspondence between spin diffeomorphism orbits of $\EDM_{(q, \lambda)}^{-1} (0)$ and $ (\EDM_{(q, \lambda)}^h)^{-1}(0)$. 

\begin{Lem}
	\label{LemLocExUn}
	The first-order prolongation of the operator $\EDM_{(q, \lambda)}^h $ is symmetric-hyperbolic (if restricted to an appropriate convex set $G$ of Lorentzian metrics), and there is a Cauchy development in the sense of \cref{DefCauchyDev}.
\end{Lem}

\begin{Prf}
	As the assertion is well known for Einstein-Maxwell theory ($\psi = 0)$, we focus on the spinor contributions only. 
	For every precompact open set $C$ we have a ball $B$ in $\Gamma(\tau^2_C)$ around $g_0|_C$ contained in $\Lor(C)$ (the Lorentzian metrics on $C$), thus there is an affine path from each metric in $B$ to $g_0$. By \cref{EqDiracghLocal} the principal symbol of $(\Dirac_g^h)^2$ can be calculated as 
	\begin{align*}
		\sigma((\Dirac^h_g)^2) (\omega,\omega) 
		= h(b_{h,g}(\omega^h), b_{h,g}(\omega^h) ) 
		= g(a_{h,g}^{-1}(\omega^h), a_{h,g}^{-1}(\omega^h)) 
		= g(\omega,\omega) 		
	\end{align*}	 
	(where, for any non-degenerate bilinear form $k$, we define $\omega^k := \omega^{\sharp_k}$).  Thus $(\Dirac_g^h)^2$ is normally hyperbolic (and consequently $\Dirac^g_h$ is symmetric-hyperbolic) as long as $g$ is globally hyperbolic.\\
	The energy momentum tensor $T_{(g,\psi)} := T_{(1, g, \psi, 0)}$ is first order in $g$, thus the second order symbol of $T_{(g,\psi)}$ w.r.t. $g$ vanishes. To calculate the first order symbol w.r.t. $\psi$, let $a$ be a smooth function with $a(p) = 0$, then $T_{(g, a \psi)} (X,Y)= 0$. The coupling of one-forms into spinors and the one from spinors into one-forms are of zeroth order, the coupling from metrics to spinors is of first order, thus one order less than the part of the operator acting on the metric. Thus all in all the highest orders vanish in the off-diagonal terms. \\
	Now the usual local existence and uniqueness holds for quasilinear symmetric-hyperbolic operators with range a convex neighborhood of $0$, \cite[16.1.2-16.1.7 and 16.2.1]{taylor3} for smooth coefficients, or \cite[Thm. 4.2]{GM} (with the same modifications as above).
\end{Prf}

The previous \cref{LemLocExUn} shows that, for every initial value set $Z= (S, g,K, \psi_0, A_0, A_1)$, the class of developments of $Z$ is nonempty.  

\begin{Lem}[\protect{geometric uniqueness, cf. \cite[Thm. 8.4]{CB}}]
\label{GeomUn}
	Let $ (M_i, g_i, \psi_i, A^{(i)}, f_i) $ be developments of the same initial values $(S, g, K, \psi_0, A_0, A_1)$. Then there are open neighborhoods $U_i$ of $f_i(S)$ in $M_i$ and a spin diffeomorphism $\tilde \alpha: \Sigma^{g_1} U_1 \to \Sigma^{g_2} U_2$ such that $ d \alpha \circ d f_1 = d f_2: T S \to TM_2$, $\alpha^* g_2 = g_1$, $\alpha^* A^{(2)} = A^{(1)}$ in Lorenz gauge and $\tilde{\alpha} \psi_1 = \psi_2$.
\end{Lem}

\begin{Prf}
	This is just an easy consequence of uniqueness of the gauged system as in \cref{LemLocExUn} and of the fact that every metric can be brought into $h$-wave gauge as in \cref{EverythingIsHarmonic}: Take $U_i$ and $f_i$ as in \cref{EverythingIsHarmonic}. By \cref{LemNaturalityIsometry} $f_i^* (g_i, A_i, \psi_i)$ are solutions of $\EDM_{(\lambda, \varepsilon)}^h=0$ to the same initial values and thus they coincide. 
\end{Prf}

\begin{Def}[Cauchy maximal]
	A Cauchy development $J = (M, g, \psi, A, f)$ is called \emf{Cauchy maximal}, if there is not proper isometric embedding $f: (M,g) \to (N,h)$ such that $f(S)$ is a Cauchy surface of $(N,h)$ if $S$ is a Cauchy surface of $(M,g)$.
\end{Def}

\begin{Thm}[existence of a Cauchy development]
	\label{ThmEDMCauchyExists}
	For a system of $k$ spinors, there is a maximal Einstein-Dirac-Maxwell Cauchy development $J = (M, g, \psi^i, A, f)$ of the initial values $Z = (S, g_0, \psi^i_0, A_0, A_1)$. If $\psi^i_0= 0$ for all $i=1, \ldots, k$ or if $A_0, A_1=0$, then $J$ is Cauchy-maximal. Furthermore, if the system is neutral and if $g$ has a conformal extension, there is a weighted Sobolev space $H$ for $\psi^i_0, A_0, A_1$ and an $\varepsilon >0$ such that if $\| ( \psi^i_0, A_0, A_1) \|_H <\varepsilon $, then $J$ is Cauchy-maximal.
\end{Thm}

\begin{Prf}
	The proof works via the usual procedure as in Choquet-Bruhat's and Geroch's seminal article \cite{CG} using \cref{GeomUn,LemLocExUn}. First, we consider the collection $C$ of all developments of $Z$. \\
	Now we need a result of Bernal and S\'anchez \cite{BS} to close a gap in the original proof given by Choquet-Bruhat and Geroch, as it was pointed out by Willie Wong \cite[p.5]{wW}. The gap is that a priori it is not clear that the collection $C$ is a set and not a proper class. Wong indicates another possible way to close the gap; nevertheless the only way manifestly without using the axiom of choice is the one via the Bernal-S\'anchez theorem, which we will need anyway for the second part of the statement.\\
	The local existence in \cref{LemLocExUn} shows that for every initial value this set of developments is nonempty. \\
	The rest is as in \cite{CG}, or as in \cite{jS, wW}, where it is shown that the axiom of choice is \emph{not} needed but only Zermelo-Fraenkel set theory axioms. \\
	Cauchy-maximality of a maximal Cauchy development $(M, g, \psi, 0)$ or $(M, g, 0, A)$ in the sense of \cref{DefCauchyDev} follows from the fact that, for any fixed globally hyperbolic metric $g$, the Dirac equation is a \emph{linear} symmetric hyperbolic system and thus has a global solution. Thus if $g$ had a Cauchy extension $(\hat{M}, \hat{g}$), then there would be a solution $\hat{\psi}$ on $(\hat{M}, \hat{g})$ to the same initial values as well, in contradiction to maximality of $(M, g, \psi, 0)$ and correspondingly for $(M,g,0,A)$. Note that here we need again the Bernal-S\'anchez result \cite{BS} on metric splittings of Lorentzian manifolds by time functions (or its refinements \cite{MS}, \cite{M}) as we do not know beforehand whether the maximal Cauchy development being a globally hyperbolic manifold is regularly sliced (and in general it is not, see \cite{M2}). The last statement follows directly from the global existence result in \cite{GM}.
\end{Prf} 

Now, one can ask the question of geodesic completeness or at least Lorentzian maximality of the maximal Cauchy development around a solution with vanishing spinor fields, e.g. around the Minkowski solution applying the same machinery as for the Klein-Gordon equation or the Maxwell equation as in \cite{LR} and \cite{jL}.

		\printbibliography[title=References]
		\addcontentsline{toc}{section}{References}
		
		\vspace{1em}
		\textsc{\textbf{Olaf Müller}, 
		Universität Regensburg, Fakultät für Mathematik, 93040 Regensburg, Germany, olaf.mueller@ur.de} \\
		
		\textsc{\textbf{Nikolai Nowaczyk}, 
		Imperial College, Maths Dept, Huxley Building, 180 Queen's Gate, Academic Visitor, Office 533, London SW7 2AZ, United Kingdom, mail@nikno.de}		
\end{document}